\documentclass[aos,nameyear,rotating,dvips]{arximspdf}
\usepackage{dcolumn}
\usepackage{multirow}
\usepackage{graphicx}

\doi{10.1214/09-AOS771}
\volume{38}
\issue{3}
\pubyear{2010}
\firstpage{1607}
\lastpage{1637}

\makeatletter

\newtheorem{prop}{Proposition}
\newproclaim{remark}{Remark}
\newproclaim{assumption}{Assumption}
\newtheorem{theorem}{Theorem}
\newtheorem{lemma}{Lemma}

\newcolumntype{d}[1]{D{.}{.}{#1}}

  \let\sv@tabnotetext\tabnotetext
  \let\sv@tabnotemark@fmt\tabnotemark@fmt
   \long\def\legend#1{{\let\tabnote@indent\leavevmode\sv@tabnotetext[]{}{#1}}}

\makeatother

\begin{document}
\begin{frontmatter}

\title{Quantile calculus and censored regression}
\runtitle{Quantile calculus and regression}

\begin{aug}
\author[A]{\fnms{Yijian} \snm{Huang}\corref{}\thanksref{t1}\ead[label=e1]{yhuang5@emory.edu}}
\runauthor{Y. Huang}
\affiliation{Emory University}
\address[A]{Department of Biostatistics\\
\quad and Bioinformatics\\
Rollins School of Public Health\\
Emory University\\
Atlanta, Georgia 30322\\
USA\\
\printead{e1}} 
\end{aug}

\thankstext{t1}{Supported in part by Grants R01 CA090747 and R01 AI078835
from the U.S. National Institutes of Health.}

\received{\smonth{8} \syear{2009}}
\revised{\smonth{11} \syear{2009}}

%
\begin{abstract}
Quantile regression has been advocated in survival analysis to assess
evolving covariate effects. However, challenges arise when the
censoring time is not always observed and may be covariate-dependent,
particularly in the presence of continuously-distributed covariates. In
spite of several recent advances, existing methods either involve
algorithmic complications or impose a probability grid. The former
leads to difficulties in the implementation and asymptotics, whereas
the latter introduces undesirable grid dependence. To resolve these
issues, we develop fundamental and general quantile calculus on
cumulative probability scale in this article, upon recognizing that
probability and time scales do not always have a one-to-one mapping
given a survival distribution. These results give rise to a novel
estimation procedure for censored quantile regression, based on
estimating integral equations. A numerically reliable and efficient
Progressive Localized Minimization (PLMIN) algorithm is proposed for
the computation. This procedure reduces exactly to the Kaplan--Meier
method in the $k$-sample problem, and to standard uncensored quantile
regression in the absence of censoring. Under regularity conditions,
the proposed quantile coefficient estimator is uniformly consistent and
converges weakly to a Gaussian process. Simulations show good
statistical and algorithmic performance. The proposal is illustrated in
the application to a clinical study.
\end{abstract}

%
\begin{keyword}[class=AMS]
\kwd[Primary ]{62N02}
\kwd[; secondary ]{62N01}.
\end{keyword}
\begin{keyword}
\kwd{Differential equation}
\kwd{estimating integral equation}
\kwd{piecewise-linear programming}
\kwd{PLMIN algorithm}
\kwd{quantile equality fraction}
\kwd{regression quantile}
\kwd{relative quantile}
\kwd{varying-coefficient model}.
\end{keyword}

\end{frontmatter}

\section{Introduction}\label{sec1}

Quantile regression [Koenker and Bassett (\citeyear{KB78})], concerning
models for conditional quantile functions, has developed into a primary
statistical methodology to investigate functional relationship between
a response and covariates. Targeting the full spectrum of quantiles, it
provides a far more complete statistical analysis than, say, classical
linear regression. This technique has a long history in econometric
applications. More recently, quantile regression has also been
advocated for survival analysis to address evolving covariate effects
which is a common phenomenon in demographic and clinical research among
others. For instance, the aging process as well as the effects of its
determinants can be vastly different at various stages of life [cf.
Koenker and Geling (\citeyear{KG01})]. On the other hand, a clinical
intervention can rarely be expected to have a constant effect, due to
the time lag in reaching full effect and to drug resistance. The
quantile regression model allows for varying regression coefficients
and thus suits these applications well. However, a main challenge
arises from censoring.

Denote the survival time by $T$ and the censoring time by $C$. As a
result of censoring, $T$ is not directly observed but through follow-up
time $X=T\wedge C$ and censoring indicator $\Delta=I(T\leq C)$, where
$\wedge$ is the minimization operator and $I(\cdot)$ is the indicator
function. Of interest is the relationship between $T$ and a $p\times1$
covariate vector $Z$ with constant 1 as the leading component. Quantile
function is an inverse of the distribution function
$F_{Z}(t)\equiv\Pr(T\leq t | {Z})$. However, ambiguities arise in the
presence of zero-density intervals; for example, zero mortality is not
uncommon at the beginning of many clinical trials since new enrollees
are relatively healthy. To be definitive, we adopt the \textit{cadlag
inverse}, that is, the inverse function that is right-continuous
with left-hand limits. The $\tau$th conditional quantile of $T$ given
$Z$ is defined as
%
%
\begin{equation}\label{qdef}
Q_{Z}(\tau)\equiv\sup\{t\dvtx F_{Z}(t)\leq\tau\},\qquad \tau\in[0,1).
\end{equation}
The quantile regression model postulates that
%
%
\begin{equation}\label{mod1}
Q_{Z}(\tau)={Z}^\top\beta_0(\tau) \qquad\forall\tau\in[0,1),
\end{equation}
where $\beta_0(\tau)$, referred to as the quantile coefficient, is a
function of
probability $\tau$. This model is semiparametric in general, but nonparametric
in the $k$-sample problem. The interest in evolving covariate effects
necessitates the functional modeling of (\ref{mod1}), which distinguishes
itself from the modeling on a specific quantile as in, for example,
median regression;
see Section \ref{sec6} for further discussion. Note that the time scale
may be, say,
logarithm-transformed, and accordingly the supports of $T$, $C$ and $X$
are not
necessarily restricted to the nonnegative half line. When all
components of
$\beta_0(\tau)$ except the intercept are constant in $\tau$, this model
reduces to the accelerated failure time mode as studied by Buckley and
James (\citeyear{BJ79}) and Tsiatis (\citeyear{T90}) among others. In
this regard, the quantile
regression model is a varying-coefficient generalization. To provide an
interpretation, suppose for the moment that model (\ref{mod1})
holds on the logarithmic time scale. Then, the effect of a covariate, other
than the leading 1 of $Z$, is to stretch or compress the baseline
survival time
(on the original scale) with a quantile-dependent stretching or compressing
factor.

With uncensored data, Koenker and Bassett (\citeyear{KB78}) generalized
sample quantile
and proposed regression quantile as a quantile coefficient estimator
via a
convex objective function.
An adaptation of the well-known Barrodale--Roberts algorithm was later
suggested by Koenker and D'Orey (\citeyear{KD87}) for the computation.
The reliability
and efficiency of this algorithm contribute to broader acceptance of the
standard methodology.
In the presence of censoring, Powell (\citeyear{P84}, \citeyear{P86})
proposed an estimation
procedure when censoring time $C$ is always observed. His approach applies
uncensored quantile regression to $X$ as the $\tau$th quantile of $X$ turns
out to be $Q_{Z}(\tau)\wedge C={Z}^\top\beta_0(\tau)\wedge C$.

However, in most survival studies, not only is the survival time subject
to censoring but also the censoring time is unobserved for uncensored
individuals. Taking the missing-data perspective of censoring, Ying, Jung
and Wei (\citeyear{YJW95}) and Honor\'{e}, Khan and Powell (\citeyear
{HKP02}) developed different
methods but with the common consistent estimation requirement
of the censoring distribution given covariates. This amounts to either an
unconditional independence censoring mechanism, or a finite-number
limitation on covariate values, or additional censoring-time modeling to
achieve a root-$n$ convergence rate of the estimated censoring distribution.
Obviously, none of these is desirable. Ying, Jung and Wei (\citeyear{YJW95})
indicated that such a restriction may be relieved by employing smoothing techniques to
nonparametrically estimate the conditional censoring distribution. As an
alternative, Wang and Wang (\citeyear{WW09}) developed a method by
nonparametrically
estimating the conditional survival distribution via kernel smoothing.
Nevertheless, Robins and Ritov (\citeyear{RR97}) argued that these
smoothing-based
methods may not be practical with moderate sample size in the presence of
multiple continuously-distributed covariates; see also
Portnoy (\citeyear{P03}).

Our investigation focuses on the same preceding data structure, but aims
to allow for generalities on both the censoring mechanism and covariates.
Specifically, we consider conditional independence censoring mechanism:
%
%
\begin{equation}\label{mod2}
T\perp C | {Z},
\end{equation}
where $\perp$ denotes statistical independence.
This problem
was first investigated by Portnoy (\citeyear{P03}), who suggested the
pivoting method
employing the redistribution-to-the-right imputation scheme for censoring
[Efron (\citeyear{E67})]. The mass of censored observations is
recursively redistributed
to adopt standard uncensored quantile regression. However, one premise
is the
quantile monotonicity so that the ``right,'' or future, is unequivocal
in the
redistribution. In the $k$-sample case, the monotonicity holds in uncensored
sample quantile, and the method reduces to the Kaplan--Meier method,
that is,
taking an inverse of the Kaplan--Meier estimator. Unfortunately, uncensored
quantile regression in general does not respect the monotonicity,
leading to
both algorithmic and analytic difficulties with Portnoy (\citeyear
{P03}). Indeed, the
asymptotic properties of the estimator have not yet been established; see
Neocleous, Vanden Branden and Portnoy (\citeyear{NBP06}). As an alternative,
Neocleous, Vanden Branden and Portnoy (\citeyear{NBP06})
advocated a closely related grid method. Most recently,
Peng and Huang (\citeyear{PH08}) proposed a functional estimating
function upon
discovering a martingale structure, and developed a grid-based quantile
coefficient estimator. Both uniform consistency and weak convergence have
been established. As for the last two methods, the grid dependence, however,
might not be completely satisfactory.

This article makes two main contributions to this problem. First of all,
fundamental and general quantile calculus is developed on probability scale,
establishing the probability-scale dynamics with allowance for zero-density
intervals and discontinuities in a distribution.
Second, from quantile calculus a well-defined estimator and a reliable and
efficient algorithm for censored quantile regression naturally emerge
on the basis of estimating integral equations. As compared with Portnoy
(\citeyear{P03}),
Neocleous, Vanden Branden and Portnoy (\citeyear{NBP06}) and Peng and
Huang (\citeyear{PH08}), this new
approach entails neither algorithmic complications nor a probability grid.
For the rest of this article, quantile calculus is presented in Section~\ref{sec2},
and the proposed estimator and algorithm in Section \ref{sec3}. The asymptotic
properties are investigated and an inference procedure suggested in
Section \ref{sec4}. Section~\ref{sec5} presents simulation results on
statistical and algorithmic
performance, and an illustration with a clinical study. Section \ref
{sec6} concludes
with discussion. The proofs are collected in Appendices
\ref{appA}--\ref{appE}.

\section{Quantile calculus}\label{sec2}

Given a survival distribution, a one-to-many mapping from probability
to time
scale may arise from zero-density intervals; adopting the cadlag definition
of quantile function is a solution given in the \hyperref
[sec1]{Introduction}. Reciprocally, a
one-to-many mapping from time to probability scale may also arise,
resulting from distributional discontinuities. Thus, time-scale theories
including counting-process martingales cannot be applied to
probability scale, unless continuity restriction is imposed on the
distribution. In uncensored quantile regression, this issue may be
bypassed by formulating the estimation as an optimization problem. However,
such an approach may not be feasible in
censored quantile regression, which calls for the development of
quantile calculus.

\subsection{The one-sample case}\label{sec21}

Drop $Z$ from the notation in this case. By definition,
$\Pr\{T<Q(\tau)\}\leq\tau\leq\Pr\{T\leq Q(\tau)\}$. Thus,
$Q(\tau)$
does not
correspond to a unique probability
$\tau$ when $\Pr\{T=Q(\tau)\}>0$. To fill in the missing piece, we introduce
the \textit{$\tau$th quantile equality fraction}:
\[
\xi(\tau)=
\frac{\tau-\Pr\{T<Q(\tau)\}}{\Pr\{T=Q(\tau)\}},
\]
which is the fraction of the probability mass at the quantile that
brings the
cumulative probability up to $\tau$. Here and
throughout, we define $0/0\equiv0$. Elementary algebra then gives
%
%
\begin{eqnarray}\label{dyn}
&&\Pr\{T<Q(\tau)\}+\Pr\{T=Q(\tau)\}\xi(\tau)\nonumber\\[-8pt]\\[-8pt]
&&\qquad=\int_0^\tau[\Pr\{T\geq Q(\nu)\}-\Pr\{T=Q(\nu)\}\xi(\nu) ]
\frac{d\nu}{1-\nu} \qquad\forall\tau\in[0,1).\nonumber
\end{eqnarray}
This result establishes the quantile dynamics on probability scale. More
significantly, it can be readily exploited to accommodate censoring.
Denote the limit of identifiability by
$\overline{\tau}=\sup\{\tau\dvtx\Pr\{C\geq Q(\tau)\}>0\}$.
\begin{prop}\label{prop1}
Suppose that $T$ and $C$ are independent and their distributions do not have
jump points in common. Consider integral equation
%
%
\begin{eqnarray}\label{eq1}
&&E \bigl(\Delta[I\{X<q(\tau)\}+I\{X=q(\tau)\}
\eta(\tau) ] \bigr)\nonumber\\[-8pt]\\[-8pt]
&&\qquad= E\int_0^\tau[I\{X\geq q(\nu)\}-I\{X=q(\nu)\}\eta(\nu) ]
\frac{d\nu}{1-\nu} \qquad\forall\tau\in[0,1),\nonumber
\end{eqnarray}
where $q(\cdot)$ is a cadlag function and $\eta(\cdot)$ takes values in
$[0,1]$.

\begin{longlist}
\item If $q(\cdot)=Q(\cdot)$ and $\eta(\tau)=\xi(\tau)$ for all
$\tau$
such that $\Pr\{T=Q(\tau)\}>0$, then (\ref{eq1}) holds.

\item If (\ref{eq1}) holds, then
%
%
\begin{eqnarray}
\label{uq1}
q(\tau) &=& Q(\tau),\\
\label{uq2}
E [\Delta I\{X=q(\tau)\}\eta(\tau) ] &=& E [\Delta I\{X=Q(\tau
)\}\xi(\tau) ]
\end{eqnarray}
for all $\tau\in(0,\overline{\tau})$.
\end{longlist}
\end{prop}
\begin{remark}\label{remark1}
The condition that the distributions of $T$ and $C$ do not share jump points
is practically needed for the identifiability of the former
and therefore the corresponding quantile function as well. The role of
$\eta(\tau)$ is to split probability mass in the case of $\Pr\{
T=Q(\tau
)\}>0$.
Equation (\ref{eq1}), however, does not determine $\eta(\tau)$ elsewhere.
But instead of, say, setting $\eta(\tau)$ to 0 in those occasions,
keeping the
more general form would be advantageous for later developments.
\end{remark}

\subsection{Quantile coefficient dynamics}\label{sec22}

Similar to the one-sample case, we assume the following assumption.
\begin{assumption}\label{assum1}
The conditional distribution functions of $T$ and $C$ given $Z$ do
not have jump points in common for all values of $Z$.
\end{assumption}

Write $\xi_{Z}(\tau)$ as the $\tau$th quantile equality fraction
for the distribution of $T$ given~$Z$. Generalize the definition of
identifiability limit as
\[
\overline{\tau}=\sup\bigl\{\tau\dvtx
E[{Z}^{\otimes2}I\{C\geq{Z}^\top\beta_0(\tau)\}]
\mbox{ is nonsingular}\bigr\},
\]
where ${v}^{\otimes2}\equiv{vv}^\top$. The one-sample
result of Proposition \ref{prop1} can then be extended.
\begin{prop}\label{prop2}
Suppose that quantile regression model (\ref{mod1}) and censoring
mechanism (\ref{mod2}) hold along with Assumption \ref{assum1}.
Consider integral equation
%
%
\begin{eqnarray}\label{iden}
&&E \bigl({Z}\Delta[I\{X<{Z}^\top\beta(\tau)\}+
I\{X={Z}^\top\beta(\tau)\}\eta_{Z}(\tau) ] \bigr)\nonumber\\
&&\qquad= E\int_0^\tau{Z} [I\{X\geq{Z}^\top\beta(\nu)\}
-I\{X={Z}^\top\beta(\nu)\}\eta_{Z}(\nu) ]
\frac{d\nu}{1-\nu}\\
\eqntext{\forall\tau\in[0,1),}
\end{eqnarray}
where $\beta(\cdot)$ is a cadlag function and $Z$-dependent $\eta
_{Z}(\cdot)$
takes values in $[0,1]$.

\begin{longlist}
\item
If $\beta(\cdot)=\beta_0(\cdot)$ and, for any given
$Z$, $\eta_{Z}(\tau)=\xi_{Z}(\tau)$ for all $\tau$ such that
$\Pr\{T=Q_{Z}(\tau) | {Z}\}>0$, then (\ref{iden}) holds.

\item In the case that both $C$ and $Z$ are discretely distributed,
if (\ref{iden}) holds, then
%
%
\begin{eqnarray}
\label{unik1}
\beta(\tau) &=& \beta_0(\tau),\\
\label{unik2}
E [{Z}\Delta I\{X={Z}^\top\beta(\tau)\}\eta_{Z}(\tau) ]
&=& E [{Z}\Delta I\{X={Z}^\top\beta_0(\tau)\}\xi_{Z}(\tau) ]
\end{eqnarray}
for all $\tau\in(0,\overline{\tau})$.
\end{longlist}
\end{prop}
\begin{remark}\label{remark2}
The admission of $\beta_0(\cdot)$ as a solution to integral
equation (\ref{iden}) is general 
in the sense
that no restriction on the survival and
censoring distributions is imposed other than Assumption \ref{assum1}. But the
uniqueness result of $\beta_0(\cdot)$ is provided only for the case
that $C$ and $Z$ are discretely distributed. It is also established
under some other conditions in Section \ref{sec4}. However, we do not
yet have a proof for the most general case.
\end{remark}
\begin{remark}\label{remark3}
Consider the censoring-absent special case with nonsingular
$E({Z}^{\otimes2})$. Then, integral equation (\ref{iden}) reduces to
\[
{D}(\tau)=\int_0^\tau\{E{Z}-{D}(\nu)\}\frac{d\nu}{1-\nu},
\]
where ${D}(\tau)$ is the left-hand side of (\ref{iden}). This
equation has a unique and closed-form solution ${D}(\tau)=\tau E{Z}$,
or
%
%
\begin{equation}\label{nocen}
E \bigl({Z} [I\{T<{Z}^\top\beta(\tau)\}
+I\{T={Z}^\top\beta(\tau)\}\eta_{Z}(\tau) ] \bigr)=\tau E{Z}.
\end{equation}
Note that $\eta_{Z}(\tau)$ affects the left-hand side at a nonsmooth
point only, that is, when $\Pr\{T={Z}^\top\beta(\tau)\}\neq0$.
With fixed
$\eta_{Z}(\tau)$,
the left-hand side may not be smooth in~$\beta(\tau)$. Nonetheless,
thanks to $\eta_{Z}(\tau)$, it can always be smooth in $\tau$ and therefore,
the equality of (\ref{nocen}) is attainable. As far as $\beta(\tau)$ is
concerned, the equation is equivalent to the minimization problem with
$E[\{T-{Z}^\top\beta(\tau)\}^-+\tau\{T-{Z}^\top\beta(\tau)\}]$, where
$a^-\equiv-aI(a<0)$. Thus, Proposition \ref{prop2} reduces to a well-known
result in
uncensored quantile regression.
\end{remark}
\begin{remark}\label{remark4}
Quantile equality fraction $\xi_{Z}(\tau)$ is a nuisance parameter. When
$\Pr({Z}={z})=0$ for a given value $z$,
$\Pr\{T=Q_{Z}(\tau) | {Z}={z}\}$ and thus $\xi_{z}(\tau)$
are not identifiable. Nevertheless, only quantity
$E[{Z}\Delta I\{X=\break {Z}^\top\beta_0(\tau)\}\xi_{Z}(\tau)]$ as a whole
is relevant to integral equation (\ref{iden}) and it is identifiable. As
evident from Remark \ref{remark3},
the notion of $\xi_{Z}(\tau)$ might not be necessary for uncensored
quantile regression, by employing minimization. But it is instrumental for
our development of censored quantile regression.
\end{remark}
\begin{remark}\label{remark5}
Proposition \ref{prop2} is more general than the martingale result of
Peng and
Huang [(\citeyear{PH08}), equation (4)], whose validity is limited to the
circumstance
of continuous survival distribution. In that special case, the former
may reduce to the latter since the mapping between
time and probability scales becomes one to one and all terms involving
$\eta_{Z}(\cdot)$ in integral equation (\ref{iden}) may vanish.
Even so, the more general form of (\ref{iden}) is still
desirable in order to derive a natural estimating integral equation by
the plug-in principle. After all, an empirical distribution is always
discrete, that is, full of discontinuities and zero-density
intervals.
\end{remark}

\subsection{Relative quantile}\label{sec23}

To facilitate probability-scale analysis both conceptually and algebraically,
we introduce the notion of \textit{relative quantile}. Anchored at the
$\tau$th
quantile, the $\{\tau+\lambda(1-\tau)\}$th quantile coefficient for
$\lambda\in[0,1)$ is referred to as the $\lambda$th relative quantile
coefficient, written as
$\beta_0(\lambda,\tau)\equiv\beta_0\{\tau+\lambda(1-\tau)\}$.
This notion provides a convenient vehicle to study quantile coefficient
$\beta_0(\cdot)$ forward from a given
probability, similar in spirit to the concept of hazard in survival
analysis.

An integral equation for $\beta_0(\lambda,\tau)$ can be derived with given
${D}(\tau)$, the left-hand side of (\ref{iden}) at $\tau$.
By algebraic manipulation, integral equation (\ref{iden}) implies
%
%
\begin{eqnarray}\label{rrq}\quad
&&E \bigl({Z}\Delta[I\{X<{Z}^\top\beta(\lambda,\tau
)\}+
I\{X={Z}^\top\beta(\lambda,\tau)\}\eta_{Z}(\lambda,\tau) ] \bigr)
-{D}(\tau)\nonumber\\
&&\qquad= E\int_0^\lambda{Z} [I\{X\geq{Z}^\top\beta(\nu,\tau)\}
-I\{X={Z}^\top\beta(\nu,\tau)\}\eta_{Z}(\nu,\tau) ]
\frac{d\nu}{1-\nu}\\
\eqntext{\forall
\lambda\in[0,1),}
\end{eqnarray}
where $\beta(\lambda,\tau)\equiv\beta\{\tau+\lambda(1-\tau)\}$ and
$\eta_{Z}(\lambda,\tau)\equiv\eta_{Z}\{\tau+\lambda(1-\tau)\}$.
Apparently,
this is a dual equation since equation (\ref{iden}) becomes a special case
when $\tau=0$.

\section{Proposed estimator and algorithm}\label{sec3}

\subsection{Estimating integral equation}\label{sec31}

The data consist of $(X_i,\Delta_i,{Z}_i)$, $i=1,\break\ldots,n$,
as $n$ i.i.d. replicates of $(X,\Delta,{Z})$. Proposition \ref{prop2}
leads naturally
to our proposed estimation procedure based on the empirical
version of integral equation~(\ref{iden}):
%
%
\begin{eqnarray}\label{edf}
&&\sum_{i=1}^n
{Z}_i\Delta_i [I\{X_i<{Z}^\top_i\beta(\tau)\}+
I\{X_i={Z}^\top_i\beta(\tau)\}w_i(\tau) ]\nonumber\\[-8pt]\\[-8pt]
&&\qquad= \sum_{i=1}^n\int_0^\tau{Z}_i [I\{X_i\geq{Z}^\top_i\beta(\nu
)\}
-I\{X_i={Z}^\top_i\beta(\nu)\}w_i(\nu) ]\frac{d\nu}{1-\nu},
\nonumber
\end{eqnarray}
where $w_i(\cdot)$ takes values in $[0,1]$. Representing a convenient
reparameterization of $\eta_{Z}(\tau)$ in (\ref{iden}), fraction
$w_i(\tau)$ serves the purpose of splitting the empirical probability mass
associated with individual $i$ when and only when $X_i={Z}^\top_i\beta
(\tau)$.
For an uncensored individual, this ensures the
continuity of $\phi_i(\tau)\equiv I\{X_i<{Z}^\top_i\beta(\tau)\}+
I\{X_i={Z}^\top_i\beta(\tau)\}w_i(\tau)$.

We shall say that $b$ \textit{interpolates} an observation $(X,\Delta
,{Z})$ if
$X=Z^\top b$.
\begin{theorem}\label{theo1}
Suppose that $\sum_{i=1}^nZ_i^{\otimes2}$ is nonsingular. Estimating integral
equation (\ref{edf}) admits a solution $\widehat{\beta}(\cdot)$
over $\tau\in
[0,1)$ with
the following properties: \textup{(i)} $\widehat{\beta}(\cdot)$ is cadlag; and
\textup{(ii)} $\widehat{\beta}(\tau)$ interpolates at least $p$ individuals
and the covariate matrix for the interpolated set is of full rank, for each
and every $\tau\in[0,1)$.
\end{theorem}
\begin{remark}\label{remark6}
Estimating integral equation (\ref{edf}) also admits a solution in the
case of
singular $\sum_{i=1}^nZ_i^{\otimes2}$, by Theorem \ref{theo1} upon eliminating
parametrization redundancy of the quantile coefficient.
\end{remark}
\begin{remark}\label{remark7}
A subtle issue concerns the fact that identifiability limit $\overline
{\tau}$
is unknown. Empirically, $\overline{\tau}$ cannot even be determined
\textit{definitively} to exceed \textit{any} $\tau>0$, which may be easily
seen in
the one-sample case. Although it is possible to estimate $\overline
{\tau
}$, we
do not terminate the estimation of $\beta_0(\cdot)$ at such an
estimate but
rather provide $\widehat{\beta}(\tau)$ for all $\tau\in[0,1)$; the
properties of
$\widehat{\beta}(\cdot)$ would otherwise become more complicated.
Precisely speaking, $\widehat{\beta}(\tau)$ is an estimator of
$\beta_0(\tau)$ provided
$\tau<\overline{\tau}$. This strategy of separating the estimation of
$\beta_0(\tau)$ provided $\tau<\overline{\tau}$ from that of
$\overline
{\tau}$
is similar to that adopted by Peng and Huang (\citeyear{PH08}). In contrast,
Portnoy (\citeyear{P03}) and Neocleous, Vanden Branden and Portnoy
(\citeyear{NBP06}) terminated their estimation
of $\beta_0(\cdot)$ once the estimate becomes nonunique, which might partly
explain the difficulties in their interval estimation.
\end{remark}

Geometrically, $\widehat{\beta}(\tau)$ for each $\tau$
is a hyperplane, partitioning the sample into
\[
\cases{
\{i\dvtx X_i<{Z}_i^\top\widehat{\beta}(\tau)\}, &\quad below set,\vspace*{2pt}\cr
\{i\dvtx X_i={Z}_i^\top\widehat{\beta}(\tau)\}, &\quad interpolated set,\vspace*{2pt}\cr
\{i\dvtx X_i>{Z}_i^\top\widehat{\beta}(\tau)\}, &\quad above set.}
\]
Each interpolated individual on the hyperplane may be split in a ratio of
$w_i(\tau)\dvtx\{1-w_i(\tau)\}$ to be associated with the below and
above sets,
respectively. This gives rise to a sample bipartition indexed by $\tau$,
and estimating integral equation (\ref{edf}) governs its evolution.

\subsection{Structuring the computation}\label{sec32}

Estimating integral equation (\ref{edf}) may be
solved exactly with the proposed Progressive Localized Minimization (PLMIN)
algorithm. The algorithm proceeds from the 0th quantile coefficient
upward in a
progressive fashion. Due to sample discreteness, $\widehat{\beta
}(\cdot)$ is piecewise
constant. We thus conveniently decompose the computation into sequential
rounds, each involving that of a 0th relative quantile coefficient
and a potential breakpoint.

Suppose that (\ref{edf}) is solved up to $\tau-$, and thus
$\phi_i(\tau-)$ of every uncensored individual is available.
Then, by continuity $\phi_i(\tau)=\phi_i(\tau-)$ of uncensored
individual $i$ is determined; obviously $\phi_i(\tau)=0$ in the case of
$\tau=0$. Inherited from the relationship
between integral equations (\ref{iden}) and (\ref{rrq}), estimating integral
equation (\ref{edf}) is equivalent to the following equation for relative
quantile coefficient:
%
%
\begin{eqnarray}\label{edfre}\quad
&&\sum_{i=1}^n
{Z}_i\Delta_i [I\{X_i<{Z}^\top_i\beta(\lambda,\tau)\}+
I\{X_i={Z}^\top_i\beta(\lambda,\tau)\}w_i(\lambda,\tau)-\phi
_i(\tau
) ]
\nonumber\\
&&\qquad= \sum_{i=1}^n\int_0^\lambda{Z}_i
[I\{X_i\geq{Z}^\top_i\beta(\nu,\tau)\}\\
&&\qquad\quad\hspace*{45.3pt}{}
-I\{X_i={Z}^\top_i\beta(\nu,\tau)\}w_i(\nu,\tau) ]\frac{d\nu
}{1-\nu},
\nonumber
\end{eqnarray}
where $w_i(\lambda,\tau)\equiv w_i\{\tau+\lambda(1-\tau)\}$. Since
$\beta(\lambda,\tau)$ remains constant from $\lambda=0$ up to a potential
relative breakpoint, say, $\lambda_b$, ${H}=\sum_{i=1}^n{Z}_i
[I\{X_i\geq{Z}^\top_i\beta(\lambda,\tau)\}
-I\{X_i={Z}^\top_i\beta(\lambda,\tau)\}
w_i(\tau)]$ is locally constant, that is, for $\lambda\in[0,\lambda_b)$.
In the case that a censored individual becomes interpolated, adopt the
convention that its $w_i(\lambda,\tau)$ remains constant locally. Write
${L}(\lambda)$ as the left-hand side of (\ref{edfre}).
Then, estimating integral equation (\ref{edfre}) is locally equivalent to
\[
{L}(\lambda)
=\int_0^\lambda\{{H}-{L}(\nu)\}\frac{d\nu}{1-\nu},\qquad
\lambda\in[0,\lambda_b),
\]
which admits a unique solution ${L}(\lambda)=\lambda{H}$ or
equivalently,
%
%
\begin{eqnarray}\label{rd0}\quad
&&\sum_{i=1}^n{Z}_i\Delta_i
[I\{X_i<{Z}^\top_i\beta(\lambda,\tau)\}+
I\{X_i={Z}^\top_i\beta(\lambda,\tau)\}w_i(\lambda,\tau)-\phi
_i(\tau
) ]
\nonumber\\[-8pt]\\[-8pt]
&&\qquad=\lambda\sum_{i=1}^n{Z}_i
[I\{X_i\geq{Z}^\top_i\beta(\lambda,\tau)\}
-I\{X_i={Z}^\top_i\beta(\lambda,\tau)\}
w_i(\tau) ].
\nonumber
\end{eqnarray}
Write $\widehat{\beta}(\lambda,\tau)\equiv\widehat{\beta}\{\tau
+\lambda(1-\tau)\}$.
Since $\widehat{\beta}
(\cdot)$
is cadlag, $\widehat{\beta}(0,\tau)$ is the solution to the above
equation with
$\lambda\downarrow0$. Subsequently, $\lambda_b$ is a $\lambda$, typically
the supremum $\lambda$, such that the equation holds with
$\beta(\lambda,\tau)=\widehat{\beta}(0,\tau)$. Furthermore,
$w_i(\lambda
_b-,\tau)$ of every
interpolated uncensored individual will be determined. Thus, solving
equation (\ref{edf}) moves forward to $\tau+\lambda_b(1-\tau)$.
The PLMIN algorithm is so named since the computation
will be conveniently carried out via minimization.

\subsection{Computing 0th relative quantile coefficient and potential
breakpoint}\label{sec33}

With the same arguments following (\ref{nocen}),
solving (\ref{rd0}) for $\widehat{\beta}(0,\tau)$ can be
reformulated as a minimization
problem:
\begin{eqnarray*}
\widehat{\beta}(0,\tau) &=&
\lim_{\lambda\downarrow0}\mathop{\arg\min}_{b}\sum_{i=1}^n
(X_i-{Z}^\top_i{b})\\
&&\hspace*{61.6pt}{} \times[I(X_i\geq{Z}^\top_i{b})
-\lambda^{-1}\Delta_i
\{I(X_i<{Z}^\top_i{b})-\phi_i(\tau)\} ],
\end{eqnarray*}
which no longer involves $w_i(\cdot,\tau)$. Further algebraic simplification
gives
%
%
\begin{equation}\label{r}
\widehat{\beta}(0,\tau)=\mathop{\arg\min}_{b}\sum
_{i=1}^n(X_i-{Z}_i^\top{b})^+
\end{equation}
subject to
\begin{eqnarray*}
&& X_i\leq{Z}^\top_i{b} \qquad\forall i\in
\mathbb{D}_-\equiv\{j\dvtx\Delta_j=1, \phi_j(\tau)=1\},
\\
&& X_i={Z}^\top_i{b} \qquad\forall i\in
\mathbb{D}_0\equiv\{j\dvtx\Delta_j=1, \phi_j(\tau)\in(0,1)\},
\\
&& X_i\geq{Z}^\top_i{b} \qquad\forall i\in
\mathbb{D}_+\equiv\{j\dvtx\Delta_j=1, \phi_j(\tau)=0\},
\end{eqnarray*}
where $a^+\equiv aI(a>0)$. For the special case of the 0th
quantile coefficient,
%
%
\begin{equation}\label{r1}
\widehat{\beta}(0) = \mathop{\arg\min}_{b}\sum_{i=1}^n(X_i-{Z}_i^\top{b})^+
\end{equation}
subject to
\[
X_i\geq{Z}^\top_i{b}\qquad
\forall i\dvtx\Delta_i=1.
\]
The minimization of (\ref{r}) is a piecewise-linear programming problem
with convex objective function, characterized by the following lemma to
Theorem \ref{theo1}. Note that, once $\widehat{\beta}(0,\tau)$ is
determined, so is $w_i(\tau )$ of a
$\widehat{\beta}(0,\tau)$-interpolated uncensored individual by
continuity of $\phi_i(\tau)$.
\begin{lemma}\label{lemma1}
Suppose that the condition of Theorem \ref{theo1} holds and that
covariates $Z_i$, $i\in\mathbb{D}_0$, are linearly independent. There
exists a minimizer $\widehat{\beta}(0,\tau)$ for problem (\ref{r}) such
that the covariate matrix for $\widehat{\beta}(0,\tau)$-interpolated
observations is of full rank. Furthermore, there exist \textup{(i)} a $p$-member
subset $\mathbb{S}$ of $\widehat{\beta}(0,\tau )$-interpolated
observations with $\mathbb{D}_0\subset\mathbb{S}$; and \textup{(ii)} for any
$\widehat{\beta}(0,\tau)$-interpolated censored individual~$i$,
\[
w_i(\tau)\in\cases{
[0,1], &\quad if $i\in\mathbb{S}$,\cr
\{0,1\}, &\quad otherwise,}
\]
such that \textup{(iii)} ${Z}_{\mathbb{S}}=\{{Z}_i\dvtx
i\in\mathbb{S}\}$ is of full rank; and \textup{(iv)}
$\widehat{H}\equiv\sum_{i=1}^n Z_i[I\{X_i\geq{Z}^\top
_i\widehat{\beta}(0,\tau )\}
-I\{X_i={Z}^\top_i\widehat{\beta}(0,\tau)\}w_i(\tau)]$ as determined
satisfies
%
%
\begin{equation}\label{mincon}
\sum_{i\in\mathbb{S}}Z_i\Delta_i\gamma_i=\widehat{H}
\end{equation}
for some $\gamma_i$, where $\gamma_i\leq0$ for $i\in\mathbb{D}_-$ and
$\gamma_i\geq0$ for $i\in\mathbb{D}_+$.
\end{lemma}

Piecewise-linear programming can be viewed as extended linear programming,
although a $\widehat{\beta}(0,\tau)$-interpolated individual may be
a censored
one and
thus not involved in the constraints. We devise an algorithm aiming
at the determination of the $p$-member interpolated subset
$\mathbb{S}$,
the same strategy as the simplex method of linear programming [e.g.,
Gill, Murray and Wright (\citeyear{GMW91})]. To locate a candidate
member of $\mathbb{S}$, the
method of steepest descent is used. Note that a feasible value for
$\widehat{\beta}(0,\tau)$ is readily available. In the case of $\tau
=0$, any value
with a
sufficiently small intercept component is feasible. Subsequently,
$\widehat{\beta}
(\tau-)$
is feasible as necessary by continuity of $\phi_i(\cdot)$ for uncensored
individuals. The minimization along a given feasible direction is
reached once an
uncensored observation becomes interpolated, or potentially so if the
interpolated observation is a censored one instead. The constrained
space is of
dimension $p$ minus the size of $\mathbb{D}_0$. For $\widehat{\beta
}(0)$, there
is no
equality constraint and the dimension is $p$. For following 0th relative
quantile coefficients, typically the dimension is 1 in which case
the minimization involves only a line search. To deal with the possibility
of more than $p$ interpolated individuals, the perturbation anti-cycling
technique in linear programming [e.g., Gill, Murray and Wright
(\citeyear{GMW91}),
Section 8.3.3] can be
adapted. In the perturbation, one may follow a tie-breaking convention to
let individuals in $\mathbb{D}_+$ precede censored ones, which in turn precede
those in $\mathbb{D}_-$. This minimization is numerically reliable and
efficient.

The minimization determines $\widehat{\beta}(0,\tau)$, $\mathbb
{S}$, $w_i(\tau)$
for each
in $\mathbb{S}$, and $\gamma_i$ for each uncensored in
$\mathbb{S}$. Plugging them into (\ref{rd0}) yields
%
%
\begin{equation}\label{brk}
\sum_{i\in\mathbb{S}}Z_i\Delta_i\{w_i(\lambda,\tau)-w_i(\tau)\}
=\lambda\sum_{i\in\mathbb{S}}Z_i\Delta_i\gamma_i.
\end{equation}
Simple algebra then gives the potential relative breakpoint
%
%
\begin{equation}\label{bp}
\lambda_b = \cases{
\displaystyle\min_{i\in\mathbb{S}\dvtx\Delta_i=1,\gamma_i\neq0}
\frac{I(\gamma_i>0)-w_i(\tau)}{\gamma_i},
&\quad $\{i\in\mathbb{S}\dvtx\Delta_i=1,\gamma_i\neq0\}\neq\varnothing
$,\cr
1, &\quad otherwise,}\hspace*{-32pt}
\end{equation}
which is proper in the sense of $0<\lambda_b\leq1$. The lower bound of
$\lambda_b$ is obvious, whereas the upper bound can easily be established
by analyzing
the intercept component of (\ref{mincon}) and (\ref{brk}). If
$\lambda_b=1$, the final quantile is reached. Otherwise, for those
uncensored,
\[
w_i(\lambda_b-,\tau) = w_i(\tau)
+\lambda_b\gamma_i,\qquad
i\in\mathbb{S}\dvtx\Delta_i=1.
\]
At least one $w_i(\lambda_b-,\tau)$ above reaches 0 or 1, so is the
corresponding $\phi_i\{\tau+\lambda_b(1-\tau)\}$. Note that
$\lambda_b$ is
a breakpoint if $\widehat{\beta}(\tau)$ interpolates exactly $p$
individuals; but not
necessarily so otherwise. Nevertheless, of importance in both cases is that
the solution moves forward in a sensible fashion.

When $\tau$ is small, $\mathbb{S}$ typically consists of uncensored
individuals
only. But as $\tau$ becomes larger, interpolated censored individuals
could emerge when $\widehat{\beta}(\tau)$ might still be uniquely
determined nonetheless.
Eventually, the computation could reach a point beyond which $\widehat
{\beta}(\tau
)$ is
no longer unique. Apparently, this phenomenon relates to the identifiability
issue; see Remark \ref{remark7}. On a different note, just like
uncensored quantile
regression, this censored quantile regression may not respect
quantile monotonicity in general.

\subsection{Relationships with standard methods in special cases}\label{sec34}

In the absence of censoring, estimating integral equation (\ref{edf}) reduces
to
%
%
\begin{equation}\label{ncee}
\sum_{i=1}^n{Z}_i [I\{T_i<{Z}_i^\top\beta(\tau)\}
+I\{T_i={Z}_i^\top\beta(\tau)\}w_i(\tau) ]=\tau\sum_{i=1}^n{Z}_i
\end{equation}
by the same approach to obtaining (\ref{nocen}) from (\ref{iden}). Thus,
$\widehat{\beta}(\cdot)$ is the cadlag function $\beta(\cdot)$ that
minimizes $\sum_{i=1}^n[\{T_i-{Z}_i^\top\beta(\tau)\}^-
+\tau\{T_i-{Z}_i^\top\beta(\tau)\}]$, reducing to
one regression quantile of Koenker and Bassett (\citeyear{KB78}); note
that the
Koenker--Bassett estimator is not always uniquely defined. In the mean time,
$1-\phi_i(\tau)$ becomes
$I\{T_i\geq{Z}_i^\top\beta(\tau)\}-I\{T_i={Z}_i^\top\beta(\tau)\}
w_i(\tau)$, which
is the regression rank score of Gutenbrunner and Jure\u{c}kov\'{a}
(\citeyear{GJ92}).

On the other\vspace*{1pt} hand, in the one-sample case, $\widehat{\beta}(\cdot)$ reduces
exactly to
the cadlag inverse of the Kaplan--Meier estimator. It is clear
from (\ref{r1})
that $\widehat{\beta}(0)$ is the first failure time and from (\ref
{bp}) that the
breakpoint is the Nelson--Aalen estimate of the hazard at $\widehat
{\beta}(0)$.
Subsequently and more generally, each estimated
0th relative quantile is a failure time and the relative
breakpoint is the Nelson--Aalen hazard estimate. In case that the last
observation is censored, the final estimated quantile is
defined as this last follow-up time by convention. More generally,
in the $k$-sample problem, $\widehat{\beta}(\cdot)$ is a linear
combination of cadlag
inverses of the $k$ Kaplan--Meier estimators.

\section{Asymptotic study and inference}\label{sec4}

In our developments thus far, we have kept our assumptions to minimal.
But the generality challenges large-sample developments in both
exposition and technicalities; see Section \ref{sec6} for further
discussion. In this
section, we shall
focus on the situation that $F_{Z}$ is continuous and free of
zero-density intervals, and additionally $C$ is continuously distributed.
These regularity conditions were also adopted in previous investigations
[Portnoy (\citeyear{P03}), Neocleous, Vanden Branden and Portnoy
(\citeyear{NBP06}),
Peng and Huang (\citeyear{PH08})]. Nevertheless,
Portnoy (\citeyear{P03}) and Neocleous, Vanden Branden and Portnoy
(\citeyear{NBP06}) required
the absence of censoring prior to and around the 0th quantile. On the other
hand, Peng and Huang (\citeyear{PH08}) presumed that the 0th quantile
is $-\infty
$. In
contrast, we do not impose any conditions on the 0th quantile.

A parameter space needs to be specified. In light of the interpolation
property of the estimator by Theorem \ref{theo1}, we require that any
$b$ in such a
parameter space satisfies
that $E[{Z}^{\otimes2}I\{{Z}^\top\beta_0(0)\leq{Z}^\top{b}\leq
{Z}^\top\beta_0(1-)\wedge C\}]$ is nonsingular. Write
eigmin as the minimum eigenvalue of a matrix. Specifically, a
parameter space containing $\beta_0(\tau)$ for all $\tau\in[0,u]$ is
given by
\begin{eqnarray*}
\mathbb{B}(u) &=&
\bigl\{{b}\in\mathbb{R}\times\mathbb{C}_{p-1}\dvtx\operatorname{eigmin}
E [{Z}^{\otimes2}I\{{Z}^\top\beta_0(0)\leq{Z}^\top{b}\\
&&\hspace*{147.5pt}\leq
{Z}^\top\beta_0(1-)\wedge C\} ]>c(u)
\bigr\},
\end{eqnarray*}
where constant $u<\overline{\tau}$, compact space
$\mathbb{C}_{p-1}\subset\mathbb{R}^{p-1}$, and positive constant
$c(u)<\operatorname{eigmin}E[{Z}^{\otimes2}I\{C\geq{Z}^\top\beta_0(u)\}]$.
Thus, all slope components are bounded but the intercept may be
$-\infty$.

Write \mbox{$\|\cdot\|$} as the Euclidean norm. Let
$\overline{F}_{Z}(t)=1-F_{Z}(t)$
and $\overline{G}_{Z}(t)=1-G_{Z}(t)=\Pr(C>t | {Z})$. Adopt the
following regularity conditions:
\begin{enumerate}[C4.]
\item[C1.] $\overline{\tau}>0$ and $\|{Z}\|$ is bounded;
\item[C2.] $F_{Z}$ and $G_{Z}$ have density functions $f_{Z}$
and $g_{Z}$, which both are continuous and bounded, uniformly
in $t$ and $Z$;
\item[C3.] $\beta_0(\cdot)$ is Lipschitz-continuous on $[\tau
_1,\tau
_2]$ for
any $\tau_1$ and $\tau_2$ such that $0<\tau_1<\tau_2<1$;
\item[C4.] there exist $u\in(0,\overline{\tau})$ and a parameter space
$\mathbb{B}(u)$ such that the maximum singular value of
\[
\Psi(b)=E \{{Z}^{\otimes2}\overline{F}_{Z}({Z}^\top{b})
g_{Z}({Z}^\top{b}) \}
[E \{{Z}^{\otimes2}\overline{G}_{Z}({Z}^\top{b})
f_{Z}({Z}^\top{b}) \} ]^{-1}
\]
is bounded uniformly in ${b}\in\mathbb{B}(u)\setminus\partial
\mathbb{B}(u)$,
where $\partial$ denotes the boundary.
\end{enumerate}
The first two conditions are self-explanatory. Conditions C3 implies
that the
survival distribution does not have zero-density intervals between $Q_{Z}(0)$
and $Q_{Z}(1-)$. Imposing constraints on censoring, condition C4 is
a sufficient and technical one to accommodate the possibility of
unbounded $\beta_0(0)$.

Throughout this section, $\widehat{\beta}(\cdot)$ is \textit{any}
cadlag solution to
estimating integral equation (\ref{edf}). The solution may not be unique,
nor is the interpolation property in Theorem \ref{theo1} necessary.
\begin{theorem}\label{theo2}
Suppose that quantile regression model (\ref{mod1}) and censoring
mechanism (\ref{mod2}) hold along with conditions \textup{C1--C4}.
Equation (\ref{iden}) implies $\beta(\tau)=\beta_0(\tau)$ for all
$\tau\in(0,u]$.
For any $l\in(0,u)$,
${\sup_{\tau\in[l,u]}}\|\widehat{\beta}(\tau)-\beta_0(\tau)\|
\rightarrow0$
almost surely.
Furthermore, $n^{1/2}\{\widehat{\beta}(\tau)-\beta_0(\tau)\}$
converges weakly to a Gaussian process on $[l,u]$.
\end{theorem}
\begin{remark}\label{remark8}
Integral equation (\ref{iden}) is an initial value problem, and estimating
integral equation (\ref{edf}) is its empirical version. Accordingly,
the large-sample study as provided in Appendix \ref{appD} exploits classical
differential equation theory and modern empirical process theory. Our
study bears similarities with that of Peng and Huang (\citeyear
{PH08}). Indeed,
under the continuity condition of C2, (\ref{edf}) is essentially
equivalent to the estimating function of Peng and Huang [(\citeyear{PH08}),
equation (5)] since $w_i(\cdot)$ becomes negligible; see also Remark \ref{remark5}.
Nevertheless, we
spare the inductive arguments of Peng and Huang (\citeyear{PH08}) in
their asymptotic\vspace*{1pt}
study as typically necessary for a grid method,
by virtue of the fact that $\widehat{\beta}(\cdot)$ is an exact
solution to~(\ref{edf}). Equally noteworthy is that the generality here
on the 0th quantile requires a more delicate treatment.
\end{remark}
\begin{remark}\label{eqpen}
In the case that ${Z}^\top\beta_0(0)$ is $-\infty$ for all $Z$, our
estimator $\widehat{\beta}(\cdot)$ is
asymptotically equivalent to that of Peng and Huang (\citeyear{PH08}) provided
that mesh
size of the grid as required by the latter is of order $o(n^{-1/2})$.
\end{remark}

To make inference, the distribution of $n^{1/2}\{\widehat{\beta
}(\cdot)-\beta
_0(\cdot
)\}$
needs to be estimated. For their estimator, Peng and Huang (\citeyear{PH08})
adapted the
resampling approach of Jin, Ying and Wei (\citeyear{JYW01}). We adopt
the same
approach by perturbing estimating integral equation (\ref{edf}). This
procedure is equivalent to a multiplier bootstrap as described in
Kosorok [(\citeyear{K08}), Section 2.2.3].
\begin{theorem}\label{theo3}
Suppose that the conditions of Theorem \ref{theo2} hold, and that
nonnegative random
variable $\xi$ has unit mean and unit variance and satisfies
$\int_0^\infty\Pr(\xi>x)^{1/2} \,dx<\infty$. Perturb estimating integral
equation (\ref{edf}) by assigning i.i.d. random variables of the
same distribution as $\xi$ and independent of the data to individuals
in the
sample, and denote a
solution to the perturbed equation by $\widehat{\beta}^*(\cdot)$. On $[l,u]$,
$n^{1/2}\{\widehat{\beta}(\cdot)-\beta_0(\cdot)\}$ has the\vspace*{1pt} same asymptotic
distribution as $n^{1/2}\{\widehat{\beta}^*(\cdot)-\widehat{\beta}(\cdot)\}$
conditionally on the data.
\end{theorem}

The standard exponential distribution, for example, may be used to generate
these perturbing random variables. By repeatedly simulating perturbed
samples, the conditional distribution of $\widehat{\beta}^*(\cdot)$
can be obtained
as an
approximation for the distribution of $\widehat{\beta}(\cdot)$.

\section{Numerical studies}\label{sec5}

The quantile regression model is formulated in $\beta_0(\cdot)$. But
alternative
covariate-effect measures can be practically useful and were used in our
application (Section \ref{sec53}). Write
\[
\mu_0(\tau_1,\tau_2)\equiv(\tau_2-\tau_1)^{-1}\int_{\tau
_1}^{\tau_2}\beta
_0(\nu)
\,d\nu.
\]
Model (\ref{mod1}) implies
%
%
\begin{equation}\label{trim}
(\tau_2-\tau_1)^{-1}\int_{\tau_1}^{\tau_2}Q_{Z}(\nu) \,d\nu
={Z}^\top\mu_0(\tau_1,\tau_2),
\end{equation}
where the left-hand side is a trimmed mean of $T$. Therefore,
$\mu_0(\tau_1,\tau_2)$ measures \textit{trimmed mean effect}. This
measure is
versatile through the choices of $\tau_1$ and $\tau_2$. In fact,
$\beta_0(\tau)=\lim_{\nu\downarrow\tau}\mu_0(\tau,\nu)$ is a
special case.
On the other hand, $\mu_0(0,1)$ is the mean effect, that is, the
regression coefficient in the linear regression model.
Originally suggested as an average effect measure by Peng and
Huang (\citeyear{PH08}), $\mu_0(\tau_1,\tau_2)$ becomes even more
appealing in
light of
its specific interpretation as revealed. With censored data,
$\mu_0(\tau_1,\tau_2)$ is identifiable when $\tau_2\leq\overline
{\tau
}$, and
a natural estimator is given by
\[
\widehat{\mu}(\tau_1,\tau_2)=(\tau_2-\tau_1)^{-1}\int_{\tau
_1}^{\tau
_2}\widehat{\beta}(\nu)
\,d\nu.
\]
Obviously, $\widehat{\mu}(\tau_1,\tau_2)$ with $0<\tau_1<\tau
_2\leq u$ is strongly
consistent and asymptotically normal under the conditions of Theorem
\ref{theo2}. The
variance can be estimated by using the simulated distribution of
$\widehat{\beta}^*(\cdot)$.
Our numerical experience suggested that
$\widehat{\mu}(\tau_1,\tau_2)$ behaves reasonably well even when
$\tau_1$
takes 0.

\subsection{Finite-sample statistical performance}\label{sec51}

Simulations were conducted to mimic a clinical trial. On the
original time scale, the baseline survival distribution was
standard exponential and the censoring distribution was uniform on $[0,5]$.
The quantile regression model held on the logarithmic time scale, with two
nonconstant covariates: $Z_1$ was Bernoulli of probability 0.5 and $Z_2$
uniform on $[0,1]$. We considered two scenarios with the following conditional
quantile functions:
\begin{eqnarray*}
Q_{Z}(\tau) &=& \log\{-\log(1-\tau)\}+(1.25\tau\wedge
0.5)Z_1+0.5Z_2,\\
Q_{Z}(\tau) &=& \log\{-\log(1-\tau)\}+0.5Z_1+0.5Z_2.
\end{eqnarray*}
Scenario 1 above involved a ramp-up effect of $Z_1$, going from none to full
linearly with probability $\tau$ and staying constant afterwards. In contrast,
scenario 2 followed the accelerated failure time model.

%
%
\begin{table}
\caption{Statistical assessment under models with two nonconstant covariates}
\label{tab1}
\tabcolsep=0pt
\begin{tabular*}{\tablewidth}{@{\extracolsep{4in minus 4in}}ld{2.0}cccd{4.0}d{4.0}cd{3.0}ccd{2.0}cc@{}}
\hline
& \multicolumn{4}{c}{\textbf{Proposed}} & \multicolumn{3}{c}{\textbf{Portnoy,
pivoting}} & \multicolumn{3}{c}{\textbf{Portnoy, grid}}
& \multicolumn{3}{c@{}}{\textbf{Peng--Huang}}\\[-4pt]
& \multicolumn{4}{c}{\hspace*{-2pt}\hrulefill} & \multicolumn{3}{c}{\hspace*{-2pt}\hrulefill}
& \multicolumn{3}{c}{\hspace*{-2pt}\hrulefill} & \multicolumn{3}{c@{}}{\hspace*{-2pt}\hrulefill}\\
$\bolds\tau$ & \multicolumn{1}{c}{\textbf{B}} & \textbf{SD} & \textbf{SE} & \textbf{CI}
& \multicolumn{1}{c}{\textbf{B}} & \textbf{SD} & \textbf{CI}
& \multicolumn{1}{c}{\textbf{B}} & \textbf{SD} & \textbf{CI} &
\multicolumn{1}{c}{\textbf{B}} & \textbf{SD} & \textbf{CI}\\
\hline
\multicolumn{14}{c}{Scenario 1: with varying covariate effect}\\[4pt]
0.1 & 1 & 521 & 551 & 93.6 & -1 & 518 & 90.5
&-51 & 526 & 90.6 & 10 & 518 & 93.5\\
& 3 & 474 & 518 & 96.2 & 3 & 474 & 95.2
& -7 & 479 & 94.7 & 4 & 475 & 95.7\\
& 7 & 794 & 866 & 95.1 & 8 & 793 & 95.2
& 6 & 805 & 93.7 & 8 & 793 & 95.1\\
0.3 & 1 & 325 & 337 & 94.4 & 3 & 324 & 92.4
&-17 & 323 & 91.8 & 16 & 323 & 94.2\\
& 3 & 311 & 325 & 94.3 & 3 & 310 & 93.9
& -3 & 312 & 93.5 & 6 & 310 & 94.8\\
& -7 & 540 & 549 & 94.9 & -9 & 539 & 93.3
& -8 & 537 & 93.5 & -9 & 541 & 93.6\\
0.5 & -6 & 254 & 258 & 93.4 & -1 & 286 & 89.8
&-19 & 253 & 90.3 & 9 & 253 & 93.2\\
& 0 & 232 & 240 & 94.9 & -3 & 238 & 92.8
& -1 & 231 & 94.2 & -1 & 231 & 95.1\\
& 5 & 408 & 414 & 94.3 & -1 & 451 & 93.2
& 4 & 408 & 92.7 & 6 & 411 & 94.2\\
0.7 & -5 & 235 & 248 & 95.9 &-16 & 294 & 90.5
&-20 & 233 & 90.6 & 12 & 240 & 95.2\\
& 9 & 220 & 239 & 95.8 & -6 & 400 & 95.0
& 6 & 213 & 95.2 & 14 & 223 & 95.9\\
& 5 & 384 & 405 & 96.0 & 53 &1096 & 94.3
& 2 & 383 & 94.1 & 11 & 391 & 96.5\\
[4pt]
\multicolumn{14}{c}{Scenario 2: accelerated failure time model}\\
[4pt]
0.1 & 3 & 506 & 534 & 94.2 & 1 & 502 & 91.2
&-50 & 511 & 90.8 & 14 & 504 & 93.5\\
& -4 & 447 & 490 & 96.6 & -3 & 445 & 95.5
& -6 & 451 & 95.5 & -3 & 444 & 95.8\\
& 8 & 756 & 820 & 95.4 & 9 & 750 & 94.4
& 6 & 768 & 94.3 & 5 & 755 & 95.4\\
0.3 & 1 & 303 & 316 & 93.8 & 1 & 302 & 92.1
&-17 & 301 & 91.8 & 16 & 302 & 94.6\\
& 2 & 270 & 286 & 94.9 & 1 & 270 & 94.5
& 0 & 269 & 94.6 & -1 & 268 & 95.5\\
& -4 & 480 & 497 & 94.7 & -6 & 481 & 93.7
& -6 & 479 & 93.7 & -5 & 481 & 94.3\\
0.5 & -5 & 252 & 254 & 93.2 & 4 & 324 & 89.8
&-19 & 250 & 90.9 & 11 & 251 & 92.7\\
& 2 & 229 & 234 & 94.7 & -2 & 254 & 92.9
& 1 & 227 & 93.7 & 1 & 229 & 94.2\\
& 5 & 405 & 404 & 94.0 & -7 & 473 & 92.4
& 4 & 404 & 92.8 & 4 & 406 & 93.7\\
0.7 & -3 & 235 & 248 & 95.7 &101 &2395 & 90.9
&-19 & 233 & 90.4 & 13 & 240 & 95.3\\
& 7 & 221 & 239 & 95.6 & 76 &2532 & 94.1
& 5 & 214 & 94.9 & 12 & 223 & 95.7\\
& 2 & 384 & 405 & 96.0 &-192&5077 & 93.8
& 0 & 384 & 94.5 & 9 & 391 & 96.5\\
\hline
\end{tabular*}
\legend{Three rows for each $\tau$ correspond to the intercept
and two slope components of estimated $\tau$th quantile coefficient.}
\legend{B: empirical bias (\mbox{$\times$}1000);
SD: empirical standard deviation ($\times1000$); SE: empirical mean
of the standard error ($\times$1000);
CI: empirical coverage (\%) of 95\% confidence interval.}
\end{table}

The sample size was 200. Under each scenario, simulations were
conducted with
1000 iterations. For both scenarios, the censoring rate was approximately
32\%. Table \ref{tab1} reports the summary statistics for the proposed
$\tau$th
quantile coefficient estimates ranging from $\tau=0.1$ to 0.7, along with
estimates based on the pivoting method of Portnoy (\citeyear{P03}), the
grid method
of Portnoy [Neocleous, Vanden Branden and Portnoy (\citeyear{NBP06})], and
Peng and Huang (\citeyear{PH08}). The two
Portnoy's methods are available in R Quantreg package, of which the latest
version at the time of this research, 4.20, was used. The default mesh size,
0.01 in this case, was adopted
for the grid method of Portnoy, and the same mesh size for Peng and
Huang. For
point estimation, the pivoting method of Portnoy had large bias and variance
at $\tau=0.7$ under both scenarios. Other than that, these estimators all
had negligible bias and similar efficiency. But the bias of the proposed
estimator seemed smaller. These findings are not surprising since
the estimator of Peng and Huang is asymptotically equivalent to the
proposed in the settings under study; see Remark \ref{eqpen}. In addition,
similar efficiency between Peng and Huang and the grid method of
Portnoy was already observed in Peng and Huang (\citeyear{PH08}).
For interval estimation, Peng and Huang employs the same procedure
as the proposed, whereas the methods of Portnoy use bootstrap. The resampling
size was set to 200 for all these methods. The standard
error of the proposed estimator agreed with the standard deviation well.
The Wald-type 95\% confidence intervals of both
the proposed and Peng and Huang achieved reasonably accurate coverage
probability. In contrast, the bootstrap confidence intervals of
Portnoy's methods had undercoverage particularly for the intercept, a finding
consistent with that reported in Peng and Huang (\citeyear{PH08}).

These preceding stimulation settings conform to the conditions of the
asymptotic study in Section \ref{sec4}. Additional settings with distributional
discontinuities and zero-density intervals of the survival time were also
considered. One such simulation involved a setting similar to the preceding
ones but having a discontinuous baseline survival distribution:
\[
Q_{Z}(\tau) = \log\{-\log(1-\tau\vee0.4)\}+(\tau\vee0.4)Z_1+0.5Z_2,
\]
where $\vee$ denotes the maximization operator. Unfortunately, the
pivoting and
grid methods of Portnoy as implemented in R Quantreg package had numerical
difficulties and their appropriate evaluation was not permitted. Both
the estimator of Peng and Huang and the proposed had negligible bias
at $\tau=0.1$ and 0.3. However, for $\tau=0.5$ and 0.7, the absolute
median-bias of Peng and Huang reached 0.136, 0.065 and 0.041 for the
intercept and two slopes, respectively. In comparison, the absolute
median-bias of the proposed estimator was no larger than 0.026 for all
three coefficients. Here, median-bias is a more appropriate summary than
mean-bias due to the skewness of these estimators resulting
from discontinuity in the survival distribution. These results were
expected since the validity of Peng and Huang is tied to the assumption of
continuous survival distribution; see Remark \ref{remark5}.

\subsection{Algorithmic performance}\label{sec52}

The proposed method was compared\break against the pivoting method of
Portnoy, the
grid method of Portnoy, and the Peng and Huang method implemented by
Koenker (\citeyear{Koenker08}). All these existing methods are implemented
with Fortran source
code in R Quantreg package. The original implementation
of Peng and Huang in R language was inappropriate for comparison due to the
inherent slower speed of R. For the two grid methods, the default mesh size,
$0.01\wedge n^{-0.7}/2$ for sample size $n$, was adopted. The proposed method
was also implemented in R with Fortran source code. To minimize
the impact of R overhead, calling the Fortran function of a method from R
was timed. The computation
was performed on a Dell 2950 computer with 3.0 GHz Pentium Xeon X5365 CPUs.

The survival time followed the accelerated failure time model with
$p-1$ nonconstant covariates
\[
\log T=\varepsilon+\sum_{m=2}^{p}\frac{(-1)^{m-1}}{2}Z^{(m)},
\]
where $\varepsilon$ followed the extreme-value distribution, and $Z^{(m)}$,
$m=2,\ldots,p$, were independent and uniformly distributed on $[0,1]$.
The number of nonconstant covariates ranged from 1 to 8, and the sample
size from 100 to 1600. Three levels of censoring, 0\%, 25\% and 50\%,
were investigated. Unless there was no censoring, the censoring time followed
the uniform distribution between 0 and a censoring rate-determined upper
bound. Computational reliability and
efficiency of various methods for point estimation of the quantile
coefficient process were assessed with 1000 iterations, shown in
Table \ref{tab2}.

%
%
\begin{sidewaystable}
\tablewidth=\textheight
\tablewidth=\textwidth
\caption{Algorithmic evaluation of timing and reliability}
\label{tab2}
\fontsize{9}{10}\selectfont{
\begin{tabular*}{\tablewidth}{@{\extracolsep{4in minus 4in}}lccd{2.1}d{2.1}d{2.1}cd{2.1}d{2.1}d{2.1}cd{3.1}d{3.1}d{3.1}cd{3.1}d{3.1}d{3.1}@{}}
\hline
& & &
\multicolumn{15}{c@{}}{\textbf{Number of nonconstant covariates}}\\[-4pt]
& & & \multicolumn{15}{c@{}}{\hrulefill}\\
\textbf{Sample size}
& & & \multicolumn{3}{c}{\textbf{1}} && \multicolumn{3}{c}{\textbf{2}}
&& \multicolumn{3}{c}{\textbf{4}} && \multicolumn{3}{c@{}}{\textbf{8}}\\
\hline
\phantom{0}100
&Prop & T & 0.6 & 0.5 & 0.5 && 0.7 & 0.6 & 0.6
&& 1.0 & 0.9 & 0.9 && 1.6 & 1.4 & 1.5\\
&PP & R & 1.9 & 1.9 & 1.6 && 1.9 & 1.9 & 1.7
&& 2.4 & 2.4 & 2.2 && 3.6 & 3.6 & 3.2\\
& & W & 0 & 0 & 4 && 0 & 1 & 11
&& 0 & 8 & 24 && 0 & 28 & 65\\
&PG & R & 3.4 & 3.4 & 2.4 && 3.9 & 3.6 & 2.8
&& 4.9 & 4.4 & 3.4 && 5.9 & 5.4 & 4.2\\
&PHK & R &47.5 &37.7 &19.7 && 42.5 &34.8 &20.5
&&40.4 &32.4 &19.4 && 32.0 &26.6 &16.3\\
\phantom{0}200
&Prop & T & 1.4 & 1.3 & 1.3 && 1.7 & 1.7 & 1.6
&& 2.4 & 2.2 & 2.3 && 4.1 & 3.9 & 4.0\\
&PP & R & 2.3 & 2.2 & 1.7 && 2.4 & 2.3 & 2.0
&& 3.2 & 3.2 & 2.9 && \multicolumn{1}{c}{---} & 4.9 & 4.4\\
& & W & 0 & 1 & 4 && 0 & 2 & 12
&& 0 & 6 & 26 && \multicolumn{1}{c}{---} & 31 & 67\\
&PG & R & 3.0 & 2.7 & 2.0 && 3.4 & 2.9 & 2.3
&& 4.5 & 3.9 & 2.9 && \multicolumn{1}{c}{---} & 4.7 & 3.6\\
&PHK & R &43.4 &31.3 &16.9 &&37.1 &28.3 &16.9
&&35.3 &27.8 &16.4 &&27.1 &20.8 &12.7\\
\phantom{0}400
&Prop & T & 4.5 & 4.3 & 4.4 && 5.5 & 5.5 & 5.8
&& 7.7 & 7.7 & 8.0 &&13.2 &13.1 &13.7\\
&PP & R & 2.7 & 2.4 & 1.8 && 2.7 & 2.5 & 2.1
&& 3.8 & 3.6 & 3.1 && \multicolumn{1}{c}{---} & 5.6 & 5.0\\
& & W & 0 & 0 & 8 && 0 & 2 & 13
&& 0 & 6 & 27 && \multicolumn{1}{c}{---} & 32 & 67\\
&PG & R & 2.8 & 2.5 & 1.8 && 3.2 & 2.7 & 2.0
&& 4.2 & 3.6 & 2.7 && \multicolumn{1}{c}{---} & 4.9 & 3.4\\
&PHK & R &39.0 &28.4 &15.1 &&33.4 &24.5 &13.8
&&31.4 &23.1 &13.9 && 23.4 &17.6 &10.8\\
\phantom{0}800
&Prop & T &16.3 &16.4 &16.9 && 20.4 &20.7 &21.4
&&29.2 &29.2 &30.1 &&48.1 &47.6 &49.5\\
&PP & R & 2.9 & 2.5 & 1.9 && 2.9 & 2.7 & 2.2
&& 3.9 & 3.7 & 3.2 && \multicolumn{1}{c}{---} & 6.1 & 5.5\\
& & W & 0 & 1 & 7 && 0 & 2 & 18
&& 0 & 7 & 27 && \multicolumn{1}{c}{---} & 30 & 67\\
&PG & R & 3.1 & 2.6 & 1.8 && 3.4 & 2.9 & 2.1
&& 4.6 & 3.8 & 2.8 && \multicolumn{1}{c}{---} & 5.3 & 4.0\\
&PHK & R &39.1 &27.9 &14.8 &&33.0 &24.1 &13.9
&&29.5 &21.9 &13.4 &&22.6 &17.1 &10.7\\
1600
&Prop & T &63.1 &63.9 &66.0 &&80.5 &79.6 &81.9
&&112.6&109.3&113.1 &&177.6&173.3&180.0\\
&PP & R & 3.0 & 2.7 & 2.0 && 2.9 & 2.7 & 2.2
&& 4.0 & 3.8 & 3.3 && \multicolumn{1}{c}{---} & 6.7 & 6.0\\
& & W & 0 & 1 & 7 && 0 & 4 & 20
&& 0 & 10 & 28 && \multicolumn{1}{c}{---} & 33 & 62\\
&PG & R & 3.5 & 2.9 & 2.0 && 3.7 & 3.1 & 2.2
&& 4.9 & 4.0 & 3.0 && \multicolumn{1}{c}{---} & 5.5 & 4.3\\
&PHK & R &38.0 &27.5 &14.7 &&30.2 &22.5 &13.6
&&27.7 &20.9 &13.0 &&21.8 &16.3 &10.5\\
\hline
\end{tabular*}
}
\legend{Prop: the proposed estimation; PP: the pivoting method
of Portnoy (\citeyear{P03}); PG: the grid method of
Portnoy [Neocleous, Vanden Branden and Portnoy
(\citeyear{NBP06})];
PHK: Peng and Huang (\citeyear{PH08}) implemented by Koenker
(\citeyear{Koenker08});
T: CPU time (millisecond) of point estimation;
R: timing relative to the proposed estimation;
W: termination rate with warning or error (\%).
Timing for PP was based on iterations free of warning
and error.
Three columns under each combination of sample size
and number of covariates correspond to 0\%,
25\% and 50\% censoring rates.
---: unavailable due to software crash.}
\end{sidewaystable}

Both the proposed and Peng and Huang methods were reliable. However,
the pivoting and grid methods of Portnoy tended to cause frequent R session
crashes in case of no censoring and more covariates. Furthermore, in the
presence of censoring, the pivoting method of Portnoy might terminate with
warning or error messages. This rate increased with the number of covariates
and censoring rate, up to 67\%.

The computer time of the proposed method was roughly constant over different
censoring levels, given sample size and number of covariates.
Comparatively, the proposed is faster than other methods
uniformly in all scenarios considered. Specifically, the pivoting
method of
Portnoy cost 1.6 to 6.7 times as much computer time, the grid
method of Portnoy cost 1.8 to 5.9 times, and the Peng and Huang method
cost 10.5 to 47.5 times. This result is remarkable since the grid methods
involve much fewer grid points than breakpoints of the proposed method at
larger sample size, suggesting that a grid point could be much more
costly to
compute.

\subsection{Application to a clinical study}\label{sec53}

For illustration, we applied the proposed estimation procedure to the Mayo
primary biliary cirrhosis dataset [Fleming and Harrington (\citeyear
{FH91}), Appendix D].
Conducted at Mayo Clinic between 1974 and 1984, the study followed
418 patients with primary biliary cirrhosis, a rare but fatal chronic liver
disease. One question of interest was concerned with prognostic factors
associated with survival. In this analysis, we considered five baseline
measures: age, edema, log(bilirubin), log(albumin) and
log(prothrombin time). Two participants had incomplete measures
and were thus removed. Our analysis data consisted of 416
patients, with a median follow-up time of 4.74 years and a censoring
rate of
61.5\%.

%
%
\begin{figure}

\includegraphics{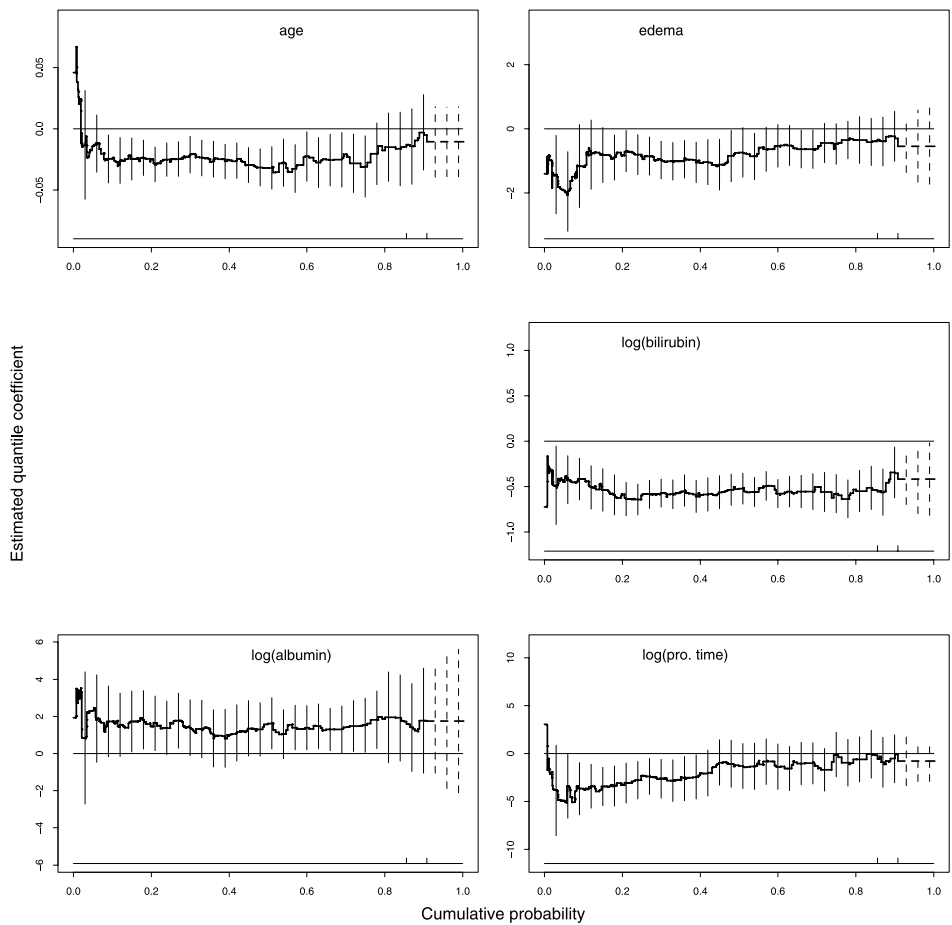}

\caption{Analysis of the Mayo primary biliary cirrhosis data.
Estimated quantile coefficient processes are shown in rugged lines, along
with pointwise Wald 95\% confidence intervals given by vertical bars.
The three regions on which the estimated quantile coefficient
hyperplanes are \textup{(a)}~unique with uncensored $\mathbb{S}$-members only,
\textup{(b)} unique
with both uncensored and censored $\mathbb{S}$-members, and
\textup{(c)} nonunique are marked on bottom horizontal lines.}
\label{fig}
\end{figure}

%
%
\begin{figure}

\includegraphics{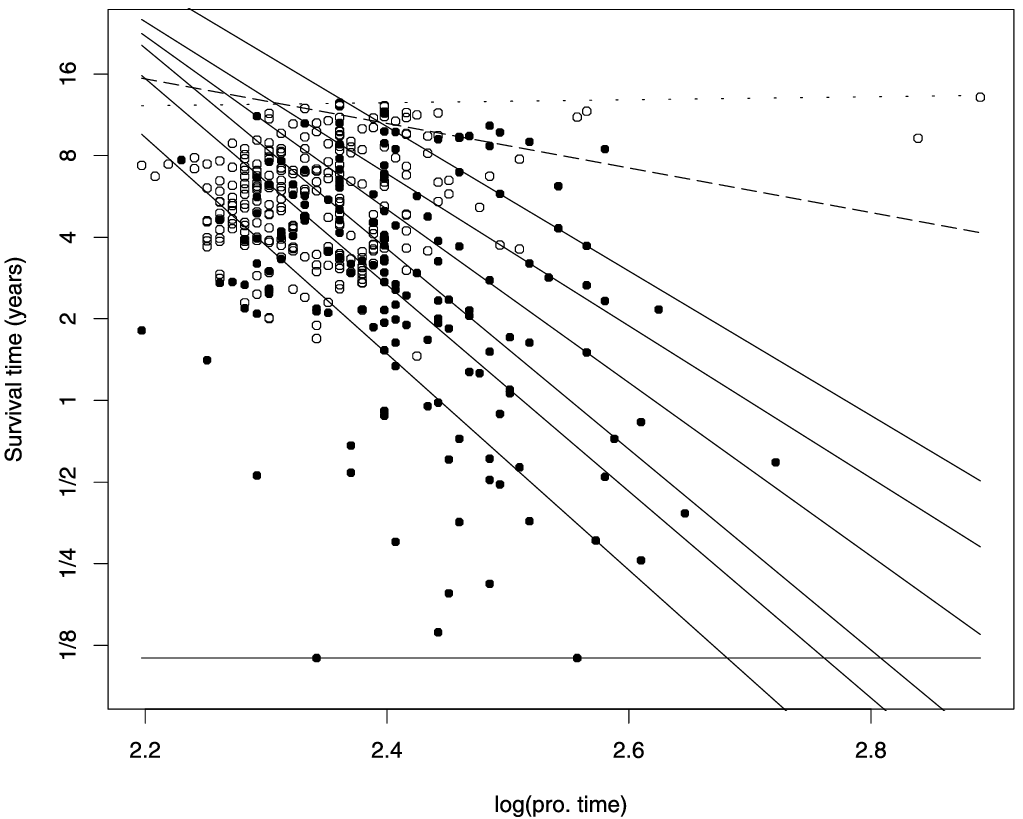}

\caption{Analysis of the Mayo primary biliary cirrhosis data with
\textup{log}(prothrombin time) as the only covariate. Dots and open circles
represent uncensored and censored individuals, respectively. Estimated
decile coefficients are shown from $\tau=0$ up to 0.8. Solid, dashed,
and dotted lines represent the corresponding hyperplanes that are
\textup{(a)}
unique with uncensored $\mathbb{S}$-members only, \textup{(b)} unique with both
uncensored and censored $\mathbb{S}$-members, and \textup{(c)} nonunique,
respectively.} \label{fig2}
\end{figure}

We adopted the quantile regression model on the logarithmic time scale, with
the five baseline measures as covariates.
The estimated quantile coefficient processes are shown
in Figure \ref{fig}, along with pointwise Wald 95\% confidence
intervals. The
resampling size for the interval estimation was 200. The maximum cumulative
probability up to which the estimated quantile coefficient was unique
is 0.91.
Among the five covariates, log(prothrombin time) in particular
exhibited a
prominent varying effect. It was negatively associated with survival
time for
short survivors, but the effect diminished gradually for
longer survivors. This result echos findings from analyses of this dataset
with other varying-coefficient models, for example, by Tian, Zucker and
Wei (\citeyear{TZW05})
using the varying-coefficient Cox model.
Nevertheless, this varying effect was not apparent in the
model with log(prothrombin time) as the only covariate, shown in
Figure \ref{fig2}.

The graphical presentation is revealing of the covariate effect evolution.
To summarize, estimated upper-trimmed mean effects and standard errors are
given in Table \ref{tab3}.
For comparison, the estimates based on the accelerated failure time model
using the log-rank and Gehan estimating functions are also included. Notice
that the two estimated regression coefficients deviate from each other for
log(prothrombin time) with the accelerated failure time model. This
disparity also suggests a lack of fit of this sub-model. In this situation,
estimates from the accelerated failure time model are difficult to
interpret. In contrast, the estimated upper-trimmed mean effects from the
quantile regression model are meaningful, for covariates with constant or
varying effects alike.

\section{Discussion}\label{sec6}

Quantile calculus as developed proves useful and effective for quantile
regression. With uncensored data, it offers a new perspective of the standard
regression procedure. Most importantly, censoring can be naturally
accommodated, and it gives rise to our proposed censored regression via
solving a well-defined estimating integral equation. To focus on the main
ideas, we have not \mbox{addressed} second-stage inference and model diagnostics,
which are practically useful. They can be developed along the lines similar
to those in Peng and Huang (\citeyear{PH08}).

%
%
\begin{table}[b]
\tabcolsep=0pt
\caption{Analysis results of the Mayo primary biliary cirrhosis data}
\label{tab3}
\begin{tabular*}{\tablewidth}{@{\extracolsep{4in minus 4in}}ld{2.4}cd{2.4}cd{2.4}cd{2.4}c@{}}
\hline
& \multicolumn{4}{c}{\textbf{Accelerated failure time model}}
& \multicolumn{4}{c@{}}{\textbf{Quantile regression model}}\\[-4pt]
& \multicolumn{4}{c}{\hspace*{-2pt}\hrulefill} & \multicolumn{4}{c@{}}{\hspace*{-2pt}\hrulefill}\\
& \multicolumn{2}{c}{\textbf{log-rank}} &
\multicolumn{2}{c}{\textbf{Gehan}}
& \multicolumn{2}{c}{$\bolds{\mu_0(0,0.8)}$}
& \multicolumn{2}{c@{}}{$\bolds{\mu_0(0,0.9)}$}\\[-4pt]
& \multicolumn{2}{c}{\hspace*{-2pt}\hrulefill} & \multicolumn{2}{c}{\hspace*{-2pt}\hrulefill}
& \multicolumn{2}{c}{\hspace*{-2pt}\hrulefill} & \multicolumn{2}{c@{}}{\hspace*{-2pt}\hrulefill}\\
& \multicolumn{1}{c}{\textbf{Est}} & \multicolumn{1}{c}{\textbf{SE}}
& \multicolumn{1}{c}{\textbf{Est}} & \multicolumn{1}{c}{\textbf{SE}}
& \multicolumn{1}{c}{\textbf{Est}} & \multicolumn{1}{c}{\textbf{SE}}
& \multicolumn{1}{c}{\textbf{Est}} & \multicolumn{1}{c@{}}{\textbf{SE}}\\
\hline
Age & -0.0259 & 0.0051 & -0.0255 & 0.0051
& -0.0238 & 0.0055 & -0.0227 & 0.0056\\
Edema & -0.7627 & 0.2276 & -0.9241 & 0.2556
& -0.8616 & 0.2413 & -0.8048 & 0.2297\\
log(bilirubin) & -0.5724 & 0.0519 & -0.5581 & 0.0611
& -0.5504 & 0.0638 & -0.5465 & 0.0615\\
log(albumin) & 1.6312 & 0.4436 & 1.4985 & 0.5013
& 1.4756 & 0.4729 & 1.4955 & 0.4438\\
log(pro.\ time) & -1.9176 & 0.5807 & -2.7761 & 0.8056
& -2.1220 & 0.8665 & -1.9426 & 0.8190\\
\hline
\end{tabular*}
\legend{Est: estimate; SE: standard error.}
\end{table}

For survival data, alternative models exist to address varying covariate
effects. One better known varying-coefficient model is the additive hazards
model of Aalen (\citeyear{A80}). There is also an extensive literature
on the
varying-coefficient Cox model, but most available estimation methods
require smoothing; see Tian, Zucker and Wei (\citeyear{TZW05}) and the
references
therein. More recently, Peng and Huang (\citeyear{PH07}) extended the
class of
semiparametric linear transformation models to allow for varying
coefficients. In comparison with all these alternatives, the quantile
regression model is particularly attractive with its simple interpretation;
see the \hyperref[sec1]{Introduction}.

When the inferential goal is on a specific quantile, the quantile regression
model for the given $\tau$ only, or the singleton model, is of direct
interest. In this case, methods for uncensored quantile regression [Koenker
and Bassett (\citeyear{KB78})] and censored one with always-observed
censoring time
[Powell (\citeyear{P84}, \citeyear{P86})] are still applicable. But
when censoring time is only
observed for censored individuals, the proposed method as well as
Portnoy (\citeyear{P03}), Neocleous, Vanden Branden and Portnoy
(\citeyear{NBP06}) and Peng and Huang (\citeyear{PH08}) may not
be applied unless the quantile regression model holds from
the 0th through the $\tau$th quantile. In contrast, the approaches of
Ying, Jung and Wei (\citeyear{YJW95}) and Honor\'{e}, Khan and Powell
(\citeyear{HKP02}) do not necessarily require
the functional model but at the price of a more restrictive censoring
mechanism. A choice between these two classes of methods depends on whether
functional survival-time modeling or censoring-time modeling might be
more reasonable and justifiable in a specific application. The method
of Wang and Wang (\citeyear{WW09}) is appealing in this regard, but
might have
practicality concerns in the case of multiple continuously-distributed
covariates, as discussed in the \hyperref[sec1]{Introduction}.\looseness=1

Generalizing the asymptotic results given in Section \ref{sec4} is of
interest, say,
to allow for zero-density intervals and discontinuities in the survival
distribution. Unfortunately, difficulties include the open question on
solution uniqueness for integral equation (\ref{iden}), as indicated in
Remark \ref{remark2}, and more. These additional ones can be readily
seen in the
one-sample case. First, the notion of consistency might not even be appropriate
in evaluating the estimated quantile corresponding to a zero-density interval.
Indeed, consistent estimation might be impossible in this circumstance but
the estimated quantile is nonetheless informative of the estimand. Second,
a~distributional discontinuity might ruin asymptotic normality of the
corresponding estimated quantile. These issues become much more
complex and also have broader implication in the general case. Due to the
sequential nature of the estimation, of concern are not only those
corresponding quantile coefficients but also the subsequent ones. Nevertheless,
it seems reasonable to conjecture that
consistency and asymptotic normality might still hold for estimated quantile
coefficients other than those corresponding to zero-density intervals and
distributional discontinuities.\looseness=1

Several additional problems are also worth further investigation. First,
our focus has been on the estimation of $\beta_0(\tau)$
provided $\tau<\overline{\tau}$, and the estimation of
identifiability limit
$\overline{\tau}$ remains to be addressed; see Remark \ref{remark7}. Second,
the submodel with a mixture of constant- and varying-coefficients would
be useful when constant effects are determined a priori for some
covariates. Efficiency gain might be expected over the more general
method in
this article. Third, in addition to right censoring, quantile
regression with
other types of censoring and truncation is also of interest.
But new techniques might be needed.
Finally, the proposed 0th
quantile coefficient estimator as given by (\ref{r1}) might be of interest
in its own right. In the absence of censoring, our estimator reduces to the
extreme regression
quantile studied by Smith (\citeyear{S94}), Portnoy and Jure\u{c}kov\'
{a} (\citeyear{PJ99})
and Chernozhukov (\citeyear{C05}) among others. Some of their results
may be extended.

\begin{appendix}
\section{\texorpdfstring{Proof of Proposition \protect\lowercase{\ref{prop1}}}{Proof of Proposition 1}}\label{appA}

Consider assertion (i). Given that $\eta(\tau)=\xi(\tau)$ for all
$\tau
$ such
that $\Pr\{T=Q(\tau)\}>0$, equation (\ref{dyn}) still holds when
$\xi
(\cdot)$
is replaced by $\eta(\cdot)$. Therefore,
\begin{eqnarray*}
&&\int_0^\tau[\Pr\{C\geq Q(\nu)\}-\Pr\{C=Q(\nu)\}\eta
(\nu
) ]\\
&&\quad{} \times d [\Pr\{T<Q(\nu)\}+\Pr\{T=Q(\nu)\}\eta(\nu
) ]\\
&&\qquad=\int_0^\tau[\Pr\{C\geq Q(\nu)\}-\Pr\{C=Q(\nu)\}\eta(\nu) ]\\
&&\qquad\quad\hspace*{12.2pt}{}
\times[\Pr\{T\geq Q(\nu)\}-\Pr\{T=Q(\nu)\}\eta(\nu) ]
\frac{d\nu}{1-\nu} \qquad\forall\tau\in[0,1).
\end{eqnarray*}
The above equation simplifies to (\ref{eq1}) under the given conditions.

For assertion (ii), only the case of $\overline{\tau}>0$ needs to be
considered. The definition of $\overline{\tau}$ implies
$E(\Delta[I\{X\geq Q(\overline{\tau})\}-I\{X=Q(\overline{\tau})\}
\xi(\overline{\tau})])=0$.
Thus,
%
%
\begin{equation}\label{a1}
E \bigl(\Delta[I\{X<Q(\overline{\tau})\}+I\{X=Q(\overline{\tau})\}
\xi(\overline{\tau}) ] \bigr)=E\Delta.
\end{equation}
Define $\tau^*=\sup\{\tau\dvtx\Pr\{C\geq q(\nu)\}>0\ \forall\nu\in
[0,\tau
]\}$.
The same argument as before gives
%
%
\begin{equation}\label{a2}
E \bigl(\Delta[I\{X<q(\tau^*)\}+I\{X=q(\tau^*)\}
\eta(\tau^*) ] \bigr)=E\Delta.
\end{equation}
Given the continuity of the left-hand side of (\ref{eq1}) in
$\tau$,
the above equation implies $\tau^*>0$. Since $\Pr\{C\geq q(\tau)\}
>0$ for
any $\tau\in[0,\tau^*)$, (\ref{eq1}) under the given conditions
implies
\begin{eqnarray*}
&&\Pr\{T<q(\tau)\}+\Pr\{T=q(\tau)\}\eta(\tau)\\
&&\qquad=\int_0^\tau[1-\Pr\{T<q(\nu)\}-\Pr\{T=q(\nu)\}\eta(\nu) ]
\frac{d\nu}{1-\nu}\qquad \forall\tau\in[0,\tau^*).
\end{eqnarray*}
The above integral equation has a unique solution:
\[
\Pr\{T<q(\tau)\}+\Pr\{T=q(\tau)\}\eta(\tau)=\tau\qquad
\forall\tau\in[0,\tau^*),
\]
from which (\ref{uq1}) and (\ref{uq2}) follow for $\tau\in(0,\tau^*)$.
Furthermore, (\ref{a1}) and (\ref{a2}) imply
$\tau^*=\overline{\tau}$. This completes the proof.

\section{\texorpdfstring{Proof of Proposition \protect\lowercase{\ref{prop2}}}{Proof of Proposition 2}}\label{appB}

Existence result (i) follows directly from Proposition \ref{prop1}. We
now prove
uniqueness result (ii) by construction. Start from $\tau=0$. Write
$\mathbb{H}$ as the discrete distributional support of $(C,{Z})$, and
define
\begin{eqnarray*}
\tau_1 &=& \sup\Bigl\{\tau\dvtx I\{c\geq{z}^\top\beta(\nu_1)\}=\lim_{\nu
_2\downarrow0}
I\{c\geq{z}^\top\beta(\nu_2)\}\mbox{ and }c\neq{z}^\top\beta
(\nu_1)\\
&&\hspace*{161.7pt}
\forall\nu_1\in(0,\tau]\mbox{ and }(c,{z})\in\mathbb{H} \Bigr\}.
\end{eqnarray*}
Thus, $I\{C\geq{Z}^\top\beta(\tau)\}$ remains constant and
$C\neq{Z}^\top\beta(\tau)$ over
$\tau\in(0,\tau_1)$ almost surely. Locally, (\ref{iden})
reduces to
\[
{D}(\tau)=\int_0^\tau\{{Y}(\tau)-{D}(\nu)\}\frac{d\nu}{1-\nu},\qquad
\tau\in[0,\tau_1),
\]
where ${D}(\tau)$ is the left-hand side of (\ref{iden}) and
${Y}(\tau)\equiv E[{Z}I\{C\geq{Z}^\top\beta(\tau)\}]$ is constant over
$\tau\in(0,\tau_1)$. The above equation admits a unique solution
${D}(\tau)=\tau{Y}(\tau)$, or equivalently
\begin{eqnarray}
E \bigl({Z}I\{C\geq{Z}^\top\beta(\tau)\} [I\{T<{Z}^\top\beta(\tau)\}
+I\{T={Z}^\top\beta(\tau)\}\eta_{Z}(\tau)-\tau] \bigr) = 0,\nonumber\\
\eqntext{\tau\in(0,\tau_1).}
\end{eqnarray}
By arguments similar to Remark \ref{remark3}, $\beta(\tau)$ is the
minimizer of
$E[\{X-{Z}^\top\beta(\tau)\wedge C\}^-+\tau\{X-{Z}^\top\beta(\tau
)\wedge C\}]$.
Recognizing that this minimization problem is the basis for
Powell's (\citeyear{P84}, \citeyear{P86}) estimator, we then obtain
(\ref{unik1}) for
$\tau\in(0,\tau_1)$ so long as $\overline{\tau}>0$.
Given (\ref{unik1}), integral equation (\ref{iden}) implies
\[
{J}(\tau)=-\int_0^\tau{J}(\nu)\frac{d\nu}{1-\nu},\qquad
\tau\in[0,\tau_1),
\]
where ${J}(\tau)$ is the difference between the two sides
of (\ref{unik2}). Thus, (\ref{unik2}) is obtained for $\tau\in
(0,\tau_1)$
by an application of the Gronwall's inequality.

Under Assumption \ref{assum1}, one can show
\begin{eqnarray*}
&&\lim_{\tau\downarrow\tau_1}E \bigl({Z} [I\{X\geq{Z}^\top
\beta
(\tau)\}
-I\{X={Z}^\top\beta(\tau)\}\eta_{Z}(\tau) ] \bigr)\\
&&\qquad=(1-\tau_1)\lim_{\tau\downarrow\tau_1}E[{Z}I\{C\geq{Z}^\top
\beta
(\tau)\}].
\end{eqnarray*}
Then, by taking advantage of the notion of relative quantile and
integral equation~(\ref{rrq}), results (\ref{unik1}) and (\ref{unik2})
can be established inductively beyond $\tau_1$, up to
$\overline{\tau}$.

\section{\texorpdfstring{Proof of Lemma \protect\lowercase{\ref{lemma1}} and
Theorem \protect\lowercase{\ref{theo1}}}{Proof of Lemma 1 and Theorem 1}}\label{appC}

With the developments in Section \ref{sec3}, it remains to establish
Lemma \ref{lemma1}.
Given the existence of a feasible value for $\widehat{\beta}(0,\tau
)$, nonnegativity
of the objective function in (\ref{r}) ensures the existence of a
minimizer. Furthermore, note that the objective function becomes linear
upon adding $X_i\leq{Z}^\top_i{b}$ or $X_i\geq{Z}^\top_i{b}$ for
each censored
individual to the constraints. Therefore, this piecewise-linear programming
problem becomes the minimization of a set of linear
programming problems, where each member involves additional constraints
concerning censored individuals. It is known that a linear programming
problem has a vertex solution if a bounded solution exists [e.g., Gill,
Murray and Wright
(\citeyear{GMW91}), Section 7.8.2]. Assertion (iii) of Lemma \ref
{lemma1} then follows.

For assertion (iv), we only consider the situation that the
interpolated set is of size $p$; otherwise one may work with the corresponding
perturbed problem.
Write $\mathbb{S}_C$ as the subset of censored individuals in $\mathbb{S}$.
The following two linear programming problems have the same solution
as (\ref{r}):
\[
\min_{b} -A^\top b
\]
subject to
\begin{eqnarray*}
X_i&\leq& Z^\top_i{b} \qquad\forall i\in\mathbb{D}_-,\qquad
X_i=Z^\top_i{b} \qquad\forall i\in\mathbb{D}_0,\\
X_i&\geq&Z^\top_i{b} \qquad\forall i\in\mathbb{D}_+,\qquad
X_i\leq Z^\top_j{b} \qquad\forall i\in\mathbb{S}_C,
\end{eqnarray*}
\[
\min_{b} -\biggl(A+\sum_{i\in\mathbb{S}_C}Z_i\biggr)^\top b
\]
subject to
\begin{eqnarray*}
X_i&\leq&Z^\top_i{b} \qquad\forall i\in\mathbb{D}_-,\qquad
X_i={Z}^\top_i{b} \qquad\forall i\in\mathbb{D}_0,\\
X_i&\geq& Z^\top_i{b} \qquad\forall i\in\mathbb{D}_+,\qquad
X_i\geq{Z}^\top_j{b} \qquad\forall i\in\mathbb{S}_C,
\end{eqnarray*}
where $A=\sum_{i\in\mathbb{D}_0}\{1-w_i(\tau)\}Z_i
+\sum_{i\in\mathbb{D}_+}Z_i+\sum_{i\dvtx\Delta_i=0,X_i>Z^\top
_i\widehat{\beta}
(0,\tau)}Z_i$.
Of course, the above two coincide in the case of $\mathbb
{S}_C=\varnothing$.
Applying Gill, Murray and Wright [(\citeyear{GMW91}), Theorem 7.8.1] yields
\[
A=\sum_{i\in\mathbb{S}}Z_i\gamma^{(1)}_i,\qquad
A+\sum_{i\in\mathbb{S}_C}Z_i=
\sum_{i\in\mathbb{S}}Z_i\gamma^{(2)}_i
\]
for some $\gamma^{(\cdot)}_i$, where $\gamma_i^{(\cdot)}\leq0$ for
$i\in\mathbb{D}_-$, $\gamma_i^{(\cdot)}\geq0$ for
$i\in\mathbb{D}_+$, and $\gamma^{(1)}_i\leq0$ and $\gamma
^{(2)}_i\geq0$
for $i\in\mathbb{S}_C$. Since $Z_{\mathbb{S}}$ is of full rank,
$\gamma_i^{(1)}=\gamma_i^{(2)}$ for $i\in\mathbb{S}\setminus
\mathbb{S}_C$
and $\gamma_i^{(1)}=\gamma_i^{(2)}-1$ for $i\in\mathbb{S}_C$.
Therefore, $\widehat{H}$ as determined upon setting
$w_i(\tau)=\gamma_i^{(2)}\in[0,1]$ for $i\in\mathbb{S}_C$ satisfies
(\ref{mincon}), with $\gamma_i=\gamma_i^{(\cdot)}$ for
$i\in\mathbb{S}\setminus\mathbb{S}_C$. This completes the proof.

\section{\texorpdfstring{Proof of Theorem \protect\lowercase{\ref{theo2}}}{Proof of Theorem 2}}\label{appD}

Similar to Peng and Huang (\citeyear{PH08}), we introduce monotone maps
$\Gamma_1({b})=E\{Z\Delta I(X\leq{Z}^\top{b})\}$ and\vspace*{1pt}
$\Gamma_2({b})=E\{Z I(X\geq{Z}^\top{b})\}$.
Write their empirical\break counterparts as
$\widehat{\Gamma}_1(b)=n^{-1}\sum_{i=1}^n Z_i\Delta_i I(X_i\leq
{Z}_i^\top
{b})$ and
$\widehat{\Gamma}_2(b)=n^{-1}\sum_{i=1}^n Z_i\times\break I(X_i\geq{Z}_i^\top{b})$.
Under condition C3, $\Gamma_1({b})$ is a one-to-one map for ${b}\in
\mathbb{B}(u)$
and $\Gamma_1^{-1}(\cdot)$ exists. Write ${H}({a})=\Gamma_2\{\Gamma
_1^{-1}(a)\}$
and note $\partial H(a)/\partial a^\top=-\Psi\{\Gamma_1^{-1}(a)\}
-\Pi
$, where
$\Psi(\cdot)$ is defined in condition C4 and $\Pi$ is the $p\times
p$ identity
matrix.

\subsection*{Identifiability}

Write $\alpha(\tau)=\Gamma_1\{\beta(\tau)\}$, and integral
equation (\ref{iden}) can be written as
\[
\alpha(\tau)=\int_0^\tau{H}\{\alpha(\nu)\}\frac{d\nu}{1-\nu};
\]
note that terms involving $\eta_{Z}(\cdot)$ vanish under condition C2.
Condition C4 implies that ${H}({a})$ is Lipschitz-continuous. The
Picard--Lindel\"{o}f theorem then asserts the solution uniqueness, that is,
$\alpha(\cdot)=\alpha_0(\cdot)\equiv\Gamma_1\{\beta_0(\tau)\}$.
It follows
that $\beta(\tau)=\beta_0(\tau)$ for all $\tau\in(0,u]$.

\subsection*{Consistency}

It is known that $\{I(X\leq{Z}^\top{b})\dvtx b\in\mathbb{R}^p\}$ is
Donsker [e.g.,
Kosorok (\citeyear{K08}), Lemma 9.12]. Furthermore, $Z$ is bounded under
condition C1. By
permanence property of the Donsker class,
$\{Z\Delta I(X\leq{Z}^\top{b})\dvtx b\in\mathbb{R}^p\}$ is Donsker. So is
$\{Z I(X\geq{Z}^\top{b})\dvtx b\in\mathbb{R}^p\}$ by similar arguments.
Since Donsker implies Glivenko--Cantelli, almost surely
%
%
\begin{equation}\label{gc}
{\sup_{{b}\in\mathbb{R}^p}}\|\widehat{\Gamma}_j(b)-\Gamma_j({b})\|=o(1),
\qquad j=1,2.
\end{equation}

On the other hand, condition C2 implies that
$\sup_{{b}\in\mathbb{R}^p}\sum_{i=1}^n I(X_i={Z}_i^\top{b})\leq p$ almost
surely. Then, coupled with condition C1, with any $w_i\in[0,1]$, almost
surely
\begin{eqnarray*}
\sup_{{b}\in\mathbb{R}^p}
\Biggl\|n^{-1}\sum_{i=1}^n{Z}_i\Delta_i I(X_i={Z}_i^\top{b})(w_i-1)
\Biggr\| &=& O(n^{-1}),\\
\sup_{{b}\in\mathbb{R}^p}
\Biggl\|n^{-1}\sum_{i=1}^n{Z}_i I(X_i={Z}_i^\top{b})w_i \Biggr\| &=& O(n^{-1}).
\end{eqnarray*}
Therefore, almost surely
%
%
\begin{equation}\label{appcon}
\sup_{\tau\in[0,u]} \biggl\|\widehat{\Gamma}_1\{\widehat{\beta}(\tau
)\}
-\int_0^\tau\widehat{\Gamma}_2\{\widehat{\beta}(\nu)\}\frac
{d\nu}{1-\nu} \biggr\|=O(n^{-1}),
\end{equation}
since (\ref{edf}) can be written as
\begin{eqnarray*}
&&\widehat{\Gamma}_1\{\beta(\tau)\}
+n^{-1}\sum_{i=1}^n{Z}_i\Delta_i I\{X_i={Z}^\top_i\beta(\tau)\}\{
w_i(\tau
)-1\}
\\
&&\qquad=\int_0^\tau\Biggl[\widehat{\Gamma}_2\{\beta(\nu)\}-n^{-1}\sum_{i=1}^n{Z}_i
I\{X_i={Z}^\top_i\beta(\nu)\}w_i(\nu) \Biggr]\frac{d\nu}{1-\nu}
\end{eqnarray*}
and $\widehat{\beta}(\cdot)$ is a solution.

Following (\ref{gc}) and (\ref{appcon}), almost surely
\[
\sup_{\tau\in[0,u]} \biggl\|\widehat{\alpha}(\tau)-\int_0^\tau{H}\{
\widehat{\alpha}(\nu)\}
\frac{d\nu}{1-\nu} \biggr\|=o(1),
\]
where $\widehat{\alpha}(\tau)=\Gamma_1\{\widehat{\beta}(\tau)\}
$. Write $L$ as the
Lipschitz constant
of ${H}(\cdot)$. Thus, for every $\epsilon>0$ and
with sufficiently large $n$, almost surely
\begin{eqnarray*}
\|\widehat{\alpha}(\tau)-\alpha_0(\tau)\| &\leq& \int_0^\tau
\|{H}\{\widehat{\alpha}(\nu)\}-{H}\{\alpha_0(\nu)\}\|\frac{d\nu
}{1-\nu}+\epsilon
\\
&\leq& \int_0^\tau L\|\widehat{\alpha}(\nu)-\alpha_0(\nu)\|
\frac{d\nu}{1-\nu}+\epsilon,
\end{eqnarray*}
which leads to
%
%
\begin{equation}\label{abound}
\|\widehat{\alpha}(\tau)-\alpha_0(\tau)\|
\leq\epsilon(1-\tau)^{-L},\qquad \tau\in[0,u],
\end{equation}
by the Gronwall's inequality.
Therefore, $\widehat{\alpha}(\tau)$ is strongly consistent for
$\alpha_0(\tau)$
uniformly in
$\tau\in[0,u]$.

It remains to show that, for any $\epsilon>0$, there exists $\delta
>0$ such
that\break ${\sup_{\tau\in[l,u]}}
\|\alpha(\tau)-\alpha_0(\tau)\|<\delta$ implies ${\sup_{\tau\in[l,u]}}
\|\beta(\tau)-\beta_0(\tau)\|<\epsilon$. Suppose that the
assertion is
false. Thus,
for each $\delta>0$, there exists $({b},\nu)$ such that
$\|\Gamma_1({b})-\alpha_0(\nu)\|<\delta$ and $\|{b}-\beta_0(\nu)\|
>d$ for
some constant $d>0$. Then, there is a subsequence of $({b},\nu)$ that
converges to, say, $({b}_0,\nu_0)$. This means that
$\Gamma_1({b}_0)=\Gamma_1\{\beta_0(\nu_0)\}$ but ${b}_0\neq\beta
_0(\nu_0)$.
However, conditions C1 and C3 imply that $f_{Z}\{{Z}^\top\beta_0(\tau
)\}$
is bounded below away from 0 uniformly in $\tau\in[l,u]$ and $Z$.
Therefore, $\partial\Gamma_1({b})/\partial{b}^\top=
E\{{Z}^{\otimes2}\overline{G}_{Z}({Z}^\top{b})
f_{Z}({Z}^\top{b})\}$ at ${b}=\beta_0(\nu_0)$ is
positive definite, which along with the monotonicity of $\Gamma
_1(\cdot)$
gives rise to a contradiction.

\subsection*{Weak convergence}

\begin{lemma}\label{lemmaA1}
Under the conditions in Theorem \ref{theo2},
%
%
\begin{eqnarray}
\label{don1}
&&{\sup_{\tau\in[0,u]}} \|\widehat{\Gamma}_1\{\widehat
{\beta}(\tau)\}-\Gamma
_1\{\widehat{\beta}
(\tau
)\}
-\widehat{\Gamma}_1\{\beta_0(\tau)\}+\Gamma_1\{\beta_0(\tau)\} \|
\nonumber\\[-8pt]\\[-8pt]
&&\qquad= o_p(n^{-1/2}),\nonumber\\
\label{don2}
&&\sup_{\tau\in[0,u]} \biggl\|\int_0^{\tau}
[\widehat{\Gamma}_2\{\widehat{\beta}(\nu)\}-\Gamma_2\{\widehat
{\beta}(\nu)\}
-\widehat{\Gamma}_2\{\beta_0(\nu)\}+\Gamma_2\{\beta_0(\nu)\} ]
\frac{d\nu}{1-\nu} \biggr\|\nonumber\\[-8pt]\\[-8pt]
&&\qquad= o_p(n^{-1/2}).
\nonumber
\end{eqnarray}
\end{lemma}
\begin{pf*}{Proof of Lemma \ref{lemmaA1}}
Consider (\ref{don1}) first. Since
$\{Z\Delta I(X\leq{Z}^\top{b})\dvtx b\in\mathbb{R}^p\}$ is Donsker,
$n^{1/2}\{\widehat{\Gamma}_1(b)-\Gamma_1({b})\}$ converges weakly to a
tight Gaussian
process. The tightness implies that, for every $\epsilon>0$ and
$m=1,\ldots,p$,
\begin{eqnarray*}
&&\lim_{\delta\downarrow0}\limsup_{n\rightarrow
\infty}
\Pr\Bigl(
\sup_{b_1,b_2 \dvtx
\sigma[n^{1/2}\{\widehat{\Gamma}_1^{(m)}(b_1)-\widehat{\Gamma
}_1^{(m)}(b_2)\}]
<\delta}n^{1/2} \bigl|\widehat{\Gamma}_1^{(m)}(b_1)\\[-4pt]
&&\hspace*{44pt}\hspace*{153.6pt}\qquad{} -\Gamma_1^{(m)}(b_1)
-\widehat{\Gamma}_1^{(m)}(b_2)\\
&&\hspace*{44pt}\hspace*{153.6pt}\hspace*{30.3pt}\qquad{}
+\Gamma_1^{(m)}(b_2) \bigr|>\epsilon\Bigr)=0,
\end{eqnarray*}
where $\sigma(\cdot)$ denotes standard deviation and superscript $^{(m)}$
the $m$th component of a vector; see, for example, Kosorok (\citeyear
{K08}). Furthermore,
note that
\begin{eqnarray*}
&&\sigma^2 \bigl[n^{1/2}\bigl\{\widehat{\Gamma}_1^{(m)}(b_1)-\widehat
{\Gamma}_1^{(m)}(b_2)\bigr\} \bigr]\\
&&\qquad= \sigma^2 \bigl[Z^{(m)}\Delta\{I(X\leq{Z}^\top b_1)-I(X\leq{Z}^\top
b_2)\} \bigr]\\
&&\qquad\leq E\bigl\{Z^{(m)2}\Delta|I(X\leq{Z}^\top b_1)-I(X\leq{Z}^\top
b_2)|\bigr\}.
\end{eqnarray*}
Write
$\Upsilon(\beta_0,\widehat{\beta},\tau)
=E[\Delta|I(X\leq{Z}^\top b)-I\{X\leq{Z}^\top\beta_0(\tau)\}|]
|_{b=\widehat{\beta}(\tau)}$.
Given condition C1, it then suffices to show
%
%
\begin{equation}\label{part1}
\sup_{\tau\in[0,u]}\Upsilon(\beta_0,\widehat{\beta},\tau)=o(1)
\end{equation}
almost surely. Let $c_f$ be the upper bound of $f_Z(\cdot)$. Apparently,
\begin{eqnarray*}
\Upsilon(\beta_0,\widehat{\beta},\tau)
&\leq& c_f E [|{Z}^\top\{b-\beta_0(\tau)\}| ]
|_{b=\widehat{\beta}(\tau)}\\
&\leq& c_f \|\widehat{\beta}(\tau)-\beta_0(\tau) \| E\|Z\|.
\end{eqnarray*}
Following the consistency of $\widehat{\beta}(\cdot)$, for every $l>0$,
%
%
\begin{equation}\label{part11}
\sup_{\tau\in[l,u]}\Upsilon(\beta_0,\widehat{\beta},\tau)=o(1)
\end{equation}
almost surely. On the other hand,
\begin{eqnarray*}
\Upsilon(\beta_0,\widehat{\beta},\tau)
&\leq& \Gamma_1^{(1)}\{\widehat{\beta}(\tau)\}+\Gamma_1^{(1)}\{
\beta_0(\tau)\}
\\
&\leq& 2\Gamma_1^{(1)}\{\beta_0(\tau)\}
+ \bigl|\Gamma_1^{(1)}\{\widehat{\beta}(\tau)\}-\Gamma_1^{(1)}\{\beta
_0(\tau)\} \bigr|.
\end{eqnarray*}
Therefore, following the consistency of $\widehat{\alpha}(\cdot)$,
for every
$\epsilon>0$,
there exists $\tau_{\epsilon}>0$ such that
%
%
\begin{equation}\label{part12}
\sup_{\tau\in[0,\tau_{\epsilon}]}\Upsilon(\beta_0,\widehat
{\beta},\tau
)<\epsilon
\end{equation}
almost surely for sufficiently large $n$. Combining (\ref{part11})
and (\ref{part12}) gives (\ref{part1}).

Now consider (\ref{don2}). Arguments similar to the above establish that,
for every $l>0$,
%
%
\begin{equation}\label{part21}\quad
{\sup_{\tau\in[l,u]}} \|\widehat{\Gamma}_2\{\widehat{\beta}(\tau
)\}-\Gamma_2\{\widehat{\beta}
(\tau
)\}
-\widehat{\Gamma}_2\{\beta_0(\tau)\}+\Gamma_2\{\beta_0(\tau)\} \|
= o_p(n^{-1/2}).
\end{equation}
Since $\sup_{b\in\mathbb{R}^p}\|\widehat{\Gamma}_2(b)-\Gamma
_2(b)\|
=O_p(n^{-1/2})$,
for every $\epsilon>0$, there exists $\tau_{\epsilon}>0$ such that
%
%
\begin{eqnarray}\label{part22}
&&\Pr\biggl(\sup_{\tau\in[0,\tau_{\epsilon}]} \biggl\|\int_0^{\tau}n^{1/2}
[\widehat{\Gamma}_2\{\widehat{\beta}(\nu)\}-\Gamma_2\{\widehat
{\beta}(\nu)\}
\nonumber\\[-8pt]\\[-8pt]
&&\hspace*{66.1pt}\qquad{}
-\widehat{\Gamma}_2\{\beta_0(\nu)\}+\Gamma_2\{\beta_0(\nu)\} ]
\frac{d\nu}{1-\nu} \biggr\|>\epsilon
\biggr)\rightarrow0.
\nonumber
\end{eqnarray}
Then, (\ref{don2}) follows from (\ref{part21}) and (\ref{part22}).
\end{pf*}

Plugging (\ref{don1}) and (\ref{don2}) into (\ref{appcon}) yields
that, for $\tau\in[0,u]$,
\begin{eqnarray*}
&&\widehat{\alpha}(\tau)-\alpha_0(\tau)
-\int_0^\tau[H\{\widehat{\alpha}(\nu)\}-{H}\{\alpha_0(\nu)\}
]\frac{d\nu}{1-\nu}
+o_p(n^{-1/2})\\
&&\qquad=-\widehat{\Gamma}_1\{\beta_0(\tau)\}
+\int_0^\tau\widehat{\Gamma}_2\{\beta_0(\nu)\}\frac{d\nu}{1-\nu},
\end{eqnarray*}
where $o_p(\cdot)$ is uniform for $\tau\in[0,u]$. Under conditions C2
and C4,
$H\{\widehat{\alpha}(\tau)\}-{H}\{\alpha_0(\tau)\}=-[\Psi\{\beta
_0(\tau)\}+\Pi+o(1)]
[\widehat{\alpha}(\tau)-\alpha_0(\tau)]$ almost surely.
Therefore,
%
%
\begin{eqnarray}\label{ord2}
&& \widehat{\alpha}(\tau)-\alpha_0(\tau)
+\int_0^\tau[\Psi\{\beta_0(\nu)\}+\Pi][\widehat{\alpha}(\nu
)-\alpha_0(\nu
)]\frac{d\nu
}{1-\nu}
\nonumber\\
&&\quad{} +o_p\bigl(n^{-1/2}+\|\widehat{\alpha}(\cdot)-\alpha_0(\cdot)\|
\bigr)\\
&&\qquad=-\widehat{\Gamma}_1\{\beta_0(\tau)\}
+\int_0^\tau\widehat{\Gamma}_2\{\beta_0(\nu)\}\frac{d\nu}{1-\nu}.
\nonumber
\end{eqnarray}

The remaining proof is sketched since it essentially follows that of
Theorem \ref{theo2}
in Peng and Huang (\citeyear{PH08}), where more details can be found.
Note that the right-hand side of (\ref{ord2}) is a martingale and converges
weakly to a Gaussian process by the martingale central limit theorem
[e.g., Fleming and Harrington (\citeyear{FH91})]. Furthermore, (\ref
{ord2}) as
a differential equation of $\widehat{\alpha}(\cdot)-\alpha_0(\cdot
)$ can be solved
by using
product integration theory [Gill and Johansen (\citeyear{GJ90})],
establishing that
$\widehat{\alpha}(\cdot)-\alpha_0(\cdot)$ as a linear map of the
right-hand side converges
weakly to a Gaussian process. The weak convergence of $\widehat{\beta
}(\cdot)$
then follows
by the functional delta method.

\section{\texorpdfstring{Proof of Theorem \protect\lowercase{\ref{theo3}}}{Proof of Theorem 3}}\label{appE}

Throughout, a quantity based on the perturbed sample is denoted by adding
an asterisk. For example, $\widehat{\Gamma}_1^*(b)$ is the perturbed
version of
$\widehat{\Gamma}_1(b)$.

The same arguments of the consistency proof in Appendix \ref{appD} may
be used to
show the strong consistency of $\widehat{\alpha}^*(\cdot)$ for
$\alpha_0(\cdot)$ on $[0,u]$
and that of $\widehat{\beta}^*(\cdot)$ for $\beta_0(\cdot)$ on
$[l,u]$, upon establishment
of the following two results. First, by Kosorok [(\citeyear{K08}),
Theorem 10.13],
almost surely
\[
{\sup_{{b}\in\mathbb{R}^p}}\|\widehat{\Gamma}_j^*(b)-\Gamma_j({b})\|=o(1),\qquad
j=1,2.
\]
Second, the terms involving $w_i(\cdot)$ in the perturbed version of
(\ref{edf}) are negligible, by the fact that the maximum of the
$n$ i.i.d. perturbing random variables is almost surely $o(n^{1/2})$
[Owen (\citeyear{O90}), Lemma 3].

By an unconditional multiplier central limit theorem [Kosorok
(\citeyear{K08}),
Theorem~10.1 and Corollary 10.3],
$n^{1/2}\{\widehat{\Gamma}_j^*(\cdot)-\Gamma_j({\cdot})\}$, $j=1,2$,
converge weakly
to tight processes. The arguments in the proof of Lemma \ref{lemmaA1}
then can be
used to establish
\begin{eqnarray*}
&&{\sup_{\tau\in[0,u]}} \|\widehat{\Gamma}_1^*\{\widehat{\beta
}^*(\tau)\}-\Gamma_1\{
\widehat{\beta}
^*(\tau)\}
-\widehat{\Gamma}_1^*\{\beta_0(\tau)\}+\Gamma_1\{\beta_0(\tau)\}
\|
\\
&&\qquad= o_p(n^{-1/2}),
\\
&&\sup_{\tau\in[0,u]} \biggl\|\int_0^{\tau}
[\widehat{\Gamma}_2^*\{\widehat{\beta}^*(\nu)\}-\Gamma_2\{
\widehat{\beta}^*(\nu)\}
-\widehat{\Gamma}_2^*\{\beta_0(\nu)\}+\Gamma_2\{\beta_0(\nu)\} ]
\frac{d\nu}{1-\nu} \biggr\|
\\
&&\qquad= o_p(n^{-1/2}).
\end{eqnarray*}
Thus, along the lines to establish (\ref{ord2}), one obtains
%
%
\begin{eqnarray}\label{ord2s}
&&\widehat{\alpha}^*(\tau)-\alpha_0(\tau)
+\int_0^\tau[\Psi\{\beta_0(\nu)\}+\Pi][\widehat{\alpha}^*(\nu
)-\alpha_0(\nu
)]\frac{d\nu}{1-\nu}
\nonumber\\
&&\quad{} +o_p\bigl(n^{-1/2}+\|\widehat{\alpha}^*(\cdot)-\alpha_0(\cdot)\|
\bigr)\\
&&\qquad=-\widehat{\Gamma}_1^*\{\beta_0(\tau)\}
+\int_0^\tau\widehat{\Gamma}_2^*\{\beta_0(\nu)\}\frac{d\nu
}{1-\nu},
\nonumber
\end{eqnarray}
given that a perturbing random variable is almost surely $o(n^{1/2})$.

Following from (\ref{ord2}) and (\ref{ord2s}),
%
%
\begin{eqnarray}\label{ord3}\quad
&&\widehat{\alpha}^*(\tau)-\widehat{\alpha}(\tau)
+\int_0^\tau[\Psi\{\beta_0(\nu)\}+\Pi][\widehat{\alpha}^*(\nu
)-\widehat{\alpha}(\nu
)]\frac{d\nu
}{1-\nu}
\nonumber\\
&&\quad{}
+o_p\bigl(n^{-1/2}+\|\widehat{\alpha}^*(\cdot)-\widehat{\alpha}(\cdot
)\|\bigr)\\
&&\qquad=-\widehat{\Gamma}_1^*\{\beta_0(\tau)\}+\widehat{\Gamma}_1\{
\beta_0(\tau)\}
+\int_0^\tau[\widehat{\Gamma}_2^*\{\beta_0(\nu)\}-\widehat
{\Gamma}_2\{\beta
_0(\nu)\}]
\frac{d\nu}{1-\nu},
\nonumber
\end{eqnarray}
since $\|\widehat{\alpha}(\cdot)-\alpha_0(\cdot)\|=O_p(n^{-1/2})$.
Note that both
$\Delta I\{X\leq{Z}^\top\beta_0(\tau)\}$ and
$\int_0^\tau I\{X\geq{Z}^\top\beta_0(\nu)\}(1-\nu)^{-1} \,d\nu$
are monotone
in $\tau$. Therefore, $\{\Delta I\{X\leq{Z}^\top\times\beta_0(\tau)\}
-\int
_0^\tau I\{X\geq{Z}^\top\beta_0(\nu)\}(1-\nu)^{-1} \,d\nu\dvtx\tau\in
[0,u]\}$
is Donsker and so is
$\{{Z}[\Delta I\{X\leq{Z}^\top\beta_0(\tau)\}-\int_0^\tau I\{X\geq
{Z}^\top
\beta_0(\nu)\}(1-\nu)^{-1} \,d\nu]\dvtx\tau\in[0,u]\}$.
By a conditional multiplier
central limit theorem [Kosorok (\citeyear{K08}), Theorem 10.4], the
right-hand side
of (\ref{ord3})
conditionally on the data converges weakly to the same Gaussian process
as the
right-hand side of (\ref{ord2}). Then, the assertion of
Theorem \ref{theo3}
follows the arguments at the end of Appendix \ref{appD}.
\end{appendix}

\section*{Acknowledgments}
The author extends his gratitude to Limin Peng for
discussions over the course of this research and her review of an earlier
version of this paper. He also thanks the seminar participants
at the University of North Carolina at Charlotte and the reviewers
for helpful comments.

%
\printaddresses


\begin{thebibliography}{99}

\bibitem[\protect\citeauthoryear{}{1980}]{A80}
\textsc{Aalen, O. O.} (1980).
A model for nonparametric regression
analysis of counting processes.
In \textit{Mathematical Statistics and Probability Theory ({P}roc.
{S}ixth {I}nternat. {C}onf., {W}is\l a, 1978)}. (W.~Klonecki, A. Kozek and J. Rosi\'{n}ski, eds.).
\textit{Lecture Notes in Statist.}
\textbf{2} 1--25.
Springer, New York.
\MR{0577267}

\bibitem[\protect\citeauthoryear{}{1979}]{BJ79}
\textsc{Buckley, J.} and \textsc{James, I.} (1979). Linear regression with
censored data. \textit{Biometrika} \textbf{66} 429--436.

\bibitem[\protect\citeauthoryear{}{2005}]{C05}
\textsc{Chernozhukov, V.}
(2005). Extremal quantile regression.
\textit{Ann. Statist.} \textbf{33} 806--839.
\MR{2163160}

\bibitem[\protect\citeauthoryear{}{1967}]{E67}
\textsc{Efron, B.} (1967). The
two sample problem with censored data.
In \textit{Proc. Fifth Berkeley Symp. Math.
Statist. Probab.} \textbf{4} 831--853. Prentice Hall, New York.

\bibitem[\protect\citeauthoryear{}{1991}]{FH91}
\textsc{Fleming, T. R.}
and \textsc{Harrington, D. P.} (1991). \textit{Counting
Processes and Survival Analysis}. Wiley, New York.
\MR{1100924}

\bibitem[\protect\citeauthoryear{}{1991}]{GMW91}
\textsc{Gill, P. E., Murray,
W.} and \textsc{Wright, M. H.} (1991). \textit{Numerical
Linear Algebra and Optimization}. Addison-Wesley, Redwood City, CA.
\MR{1074004}


\bibitem[\protect\citeauthoryear{}{1990}]{GJ90} \textsc{Gill, R.
D.} and
\textsc{Johansen, S.} (1990). A survey of product-integration
with a view towards application in survival analysis. \textit{Ann. Statist.}
\textbf{18} 1501--1555.
\MR{1074422}

\bibitem[\protect\citeauthoryear{}{1992}]{GJ92} \textsc
{Gutenbrunner, C.}
and \textsc{Jure\u{c}kov\'{a}, J.} (1992). Regression
rank scores and regression quantiles. \textit{Ann. Statist.}
\textbf{20} 305--330.
\MR{1150346}

\bibitem[\protect\citeauthoryear{}{2002}]{HKP02} \textsc{Honor\'{e},
B., Khan, S.} and \textsc{Powell, J. L.} (2002). Quantile
regression under random censoring. \textit{J. Econometrics} \textbf{109}
67--105.
\MR{1899693}

\bibitem[\protect\citeauthoryear{}{2001}]{JYW01} \textsc{Jin, Z.,
Ying, Z.}
and \textsc{Wei, L. J.} (2001). A simple resampling
method by perturbing the minimand. \textit{Biometrika} \textbf{88} 381--390.
\MR{1844838}

\bibitem[\protect\citeauthoryear{}{2008}]{Koenker08}
\textsc{Koenker, R.} (2008).
Censored quantile regression
redux. \textit{J. Stat. Softw.} \textbf{27} 1--25.

\bibitem[\protect\citeauthoryear{}{1978}]{KB78}
\textsc{Koenker, R.} and \textsc{Bassett, G.} (1978). Regression quantiles.
\textit{Econometrica} \textbf{46} 33--50.
\MR{0474644}

\bibitem[\protect\citeauthoryear{}{1987}]{KD87}
\textsc{Koenker, R.} and \textsc{D'Orey, V.} (1987). Computing regression
quantiles. \textit{Appl. Statist.} \textbf{36} 383--393.

\bibitem[\protect\citeauthoryear{}{2001}]{KG01}
\textsc{Koenker, R.} and \textsc{Geling, O.} (2001). Reappraising medfly
longevity: A quantile regression survival analysis.
\textit{J. Amer. Statist. Assoc.} \textbf{96} 458--468.
\MR{1939348}

\bibitem[\protect\citeauthoryear{}{2008}]{K08}
\textsc{Kosorok, M. R.}
(2008). \textit{Introduction to Empirical Processes
and Semiparametric Inference}.  Springer, New York.

\bibitem[\protect\citeauthoryear{}{2006}]{NBP06} \textsc{Neocleous,
T., Vanden
Branden, K.} and \textsc{Portnoy, S.} (2006).
Correction to ``{C}ensored regression quantiles,'' by S. Portnoy.
\textit{J. Amer. Statist. Assoc.} \textbf{101} 860--861.
\MR{2281250}

\bibitem[\protect\citeauthoryear{}{1990}]{O90}
\textsc{Owen, A.} (1990).
Empirical likelihood ratio confidence
regions. \textit{Ann. Statist.} \textbf{18} 90--120.
\MR{1041387}

\bibitem[\protect\citeauthoryear{}{2007}]{PH07}
\textsc{Peng, L.} and
\textsc{Huang, Y.} (2007). Survival analysis with temporal
covariate effects. \textit{Biometrika} \textbf{94} 719--733.
\MR{2410019}

\bibitem[\protect\citeauthoryear{}{2008}]{PH08}
\textsc{Peng, L.} and \textsc{Huang, Y.} (2008). Survival analysis
with quantile
regression models. \textit{J. Amer. Statist. Assoc.} \textbf{103} 637--649.
\MR{2435468}

\bibitem[\protect\citeauthoryear{}{2003}]{P03}
\textsc{Portnoy, S.} (2003).
Censored regression quantiles.
\textit{J. Amer. Statist. Assoc.} \textbf{98} 1001--1012.
\MR{2041488}

\bibitem[\protect\citeauthoryear{}{1999}]{PJ99}
\textsc{Portnoy, S.} and \textsc{Jure\u{c}kov\'{a}, J.} (1999). On
extreme regression quantiles. \textit{Extremes} \textbf{2} 227--243.
\MR{1781938}

\bibitem[\protect\citeauthoryear{}{1984}]{P84}
\textsc{Powell, J. L.} (1984).
Least absolute deviations estimation for
the censored regression model. \textit{J.~Econometrics} \textbf{25} 303--325.
\MR{0752444}

\bibitem[\protect\citeauthoryear{}{1986}]{P86}
\textsc{Powell, J. L.} (1986).
Censored regression quantiles.
\textit{J. Econometrics} \textbf{32} 143--155.
\MR{0853049}

\bibitem[\protect\citeauthoryear{}{1997}]{RR97}
\textsc{Robins, J. M.} and \textsc{Ritov, Y.} (1997). Toward a curse of
dimensionality appropriate (CODA) asymptotic theory for semi-parametric
models. \textit{Stat. Med.} \textbf{16} 285--319.

\bibitem[\protect\citeauthoryear{}{1994}]{S94}
\textsc{Smith, R.} (1994).
Nonregular regression. \textit{Biometrika}
\textbf{81} 173--183.
\MR{1279665}

\bibitem[\protect\citeauthoryear{}{2005}]{TZW05}
\textsc{Tian, L., Zucker, D.}
and \textsc{Wei, L. J.} (2005). On the Cox model
with time-varying regression coefficients. \textit{J. Amer. Statist. Assoc.}
\textbf{100} 172--183.
\MR{2156827}

\bibitem[\protect\citeauthoryear{}{1990}]{T90}
\textsc{Tsiatis, A. A.}
(1990). Estimating regression parameters using
linear rank tests for censored data. \textit{Ann. Statist.} \textbf{18}
354--372.
\MR{1041397}

\bibitem[\protect\citeauthoryear{}{2009}]{WW09}
\textsc{Wang, H. J.} and \textsc{Wang, L.} (2009). Locally weighed
censored quantile
regression. \textit{J. Amer. Statist. Assoc.} \textbf{104} 1117--1128.
\MR{2562007}

\bibitem[\protect\citeauthoryear{}{1995}]{YJW95}
\textsc{Ying, Z., Jung, S. H.} and \textsc{Wei, L. J.} (1995).
Survival analysis with median regression models.
\textit{J.~Amer. Statist. Assoc.}
\textbf{90} 178--184.
\MR{1325125}

\end{thebibliography}
\end{document}